\documentclass[11pt]{article}
\usepackage[T1]{fontenc}

\usepackage{fullpage}
\usepackage{amsthm}
\usepackage{amssymb}
\usepackage{amsmath}
\usepackage{graphicx}
\usepackage{tikz}
\usepackage{mathrsfs}
\usepackage{float}
\usepackage{makecell}
\usepackage{amscd}
\usepackage{subcaption}
\usepackage{multirow}
\usepackage{cite}
\usepackage{verbatim}
\usepackage{hhline}
\usepackage{comment}
\usepackage{paralist}

\usepackage{hyperref}

\usepackage{graphicx}
\usepackage{caption}

\usepackage{pgfplots}
\pgfplotsset{compat=1.15}
\usetikzlibrary{arrows}

\renewcommand{\epsilon}{\varepsilon}

\usepackage{relsize}
\newtheorem{theorem}{Theorem}
\newtheorem{lemma}[theorem]{Lemma}

\newtheorem{proposition}[theorem]{Proposition}
\newtheorem{corollary}[theorem]{Corollary}

\newtheorem{definition}[theorem]{Definition}

\theoremstyle{remark}
\newtheorem{remark}[theorem]{Remark}
\theoremstyle{plain}

\newcommand{\bb}{\begin{equation}}
\newcommand{\ee}{\end{equation}}

\def\IB{\relax\hbox{$\inbar\kern-.3em{\rm B}$}}
\def\IC{\relax\hbox{$\inbar\kern-.3em{\rm C}$}}
\def\ID{\relax\hbox{$\inbar\kern-.3em{\rm D}$}}
\def\IE{\relax\hbox{$\inbar\kern-.3em{\rm E}$}}
\def\IF{\relax\hbox{$\inbar\kern-.3em{\rm F}$}}
\def\IG{\relax\hbox{$\inbar\kern-.3em{\rm G}$}}
\def\IGa{\relax\hbox{${\rm I}\kern-.18em\Gamma$}}
\def\IH{\relax{\rm I\kern-.18em H}}
\def\IK{\relax{\rm I\kern-.18em K}}
\def\IL{\relax{\rm I\kern-.18em L}}
\def\IP{\relax{\rm I\kern-.18em P}}
\def\IR{\relax{\rm I\kern-.18em R}}

\def\IZ{\relax{\rm Z\kern-.5em Z}}

\title{\vspace{-2em}Counting Triangulated 3-Spheres, Nuclei, and Polyhedra}
\date{}
\author{}

\begin{document}
\author{Bergfinnur Durhuus\thanks{University of Copenhagen, Denmark} \and Marta Pavelka\thanks{University of Copenhagen, Denmark}}

\maketitle

\begin{abstract}
We exploit two-dimensional dual polyhedra to study
triangulated $3$-spheres and nuclei.
First, we show that the long-standing open problem of whether
the number of triangulated $3$-spheres is exponentially
bounded is equivalent to the corresponding problem for
one-vertex generalized triangulations, and hence also
for special spines of the $3$-ball. This also establishes
an equivalent counting problem for nuclei.
Second, we investigate properties of nuclei, including
their dual polyhedra. We also prove that every nontrivial
nucleus contains a tetrahedron with six interior edges,
identify the unique vertex-minimal nontrivial nucleus
with $11$ vertices and $35$ tetrahedra, and construct
nuclei whose dual graphs have genus growing linearly
with the number of tetrahedra.
\end{abstract}

\section{Introduction}
The problem of enumerating triangulations (or quadrangulations) of the $2$-sphere (or planar maps, more generally) up to combinatorial equivalence was originally considered and solved by W. T. Tutte
\cite{Tutte1962Census, Tutte1963Census}. Over the years, a variety of techniques, including analytic combinatorics and matrix integral approaches, have been developed, leading to detailed results on a large class of similar problems, like  the enumeration of various classes of graphs and discretized surfaces of fixed genus \cite{Schaeffer1998, Ambjorn1997QuantumGeometry}.  These results have contributed to progress in understanding the statistical mechanics of $2$-dimensional systems (see e.g. \cite{Ambjorn1997QuantumGeometry}) and the construction of novel probability distributions on continuous surfaces \cite{LeGall2013BrownianMap,HoldenSun2026SIAM} related to $2$-dimensional Liouville quantum gravity \cite{Polyakov1981,DavidKupiainenRhodesVargas2016}.

A general feature of the mentioned combinatorial results is that restricting to surfaces of fixed genus yields an exponential upper bound on the number of surfaces as a function of the number of $2$-cells, whereas it is generally easy to establish super-exponential lower bounds for surfaces of arbitrary genus. Attempts to obtain similar exponential bounds on topological classes of higher dimensional geometric objects have encountered serious obstacles, and relatively few applicable tools are available.

In particular, the question whether the number of triangulated $3$-spheres (up to combinatorial equivalence) is exponentially bounded \cite{Ambjorn1991ThreeDimensional,DJ95,Gromov2000SpacesQuestions,BZ11,BP24} is still open. Considerable effort has gone into establishing exponential bounds for restricted classes of triangulated $3$-spheres, with the hope that such classes might encompass all triangulated $3$-spheres. The first author and Jónsson proved an exponential bound for locally constructible (LC) triangulations~\cite{DJ95} and conjectured that all triangulated $3$-spheres are LC. This was disproved by Benedetti and Ziegler~\cite{BZ11}. More recently, Benedetti and the second author established an exponential bound for the larger class of $2$-LC triangulated manifolds~\cite{BP24}. Although this class does not contain all triangulated $3$-balls \cite{BB2017}, no example of a triangulated $3$-sphere that is not $2$-LC is currently known. The obstructions to the $2$-LC construction of
$3$-balls lead naturally to the class of
triangulations studied in this paper.

A main result of this paper is an equivalence between the problem of exponentially bounding the number of triangulated $3$-spheres and the corresponding problem for certain contractible $2$-dimensional polyhedra, called \emph{special spines} of the $3$-ball
\[
B^3=\{x\in\mathbb R^3:\|x\|\le1\}.
\]

Special spines have a long history in the study of $3$-manifold topology~\cite{Matveev2007}. Unlike the polyhedral surfaces familiar from combinatorics, the special spines considered here belong to a class of $2$-dimensional polyhedra called \emph{fake surfaces}, in which three half-planes may meet along edges. The combinatorial structure of fake surfaces has also attracted renewed interest. Recently, Fagan, Qiu, and Wang~\cite{FaganQiuWang2026} classified acyclic cellular fake surfaces with up to four vertices. Among the $4618$ such polyhedra with four vertices, only $128$ are special spines. They conjecture that the proportion of special spines among contractible fake surfaces tends to zero as the number of vertices increases.

A key step in establishing the mentioned equivalence is a reduction from simplicial triangulations of the $3$-sphere
\[
S^3=\{x\in\mathbb R^4:\|x\|=1\}
\]
to one-vertex generalized triangulations, using Matveev's arch construction \cite[p.~7]{Matveev2007}. The second main aspect of this paper concerns a subclass of triangulations of the 3-ball, called nuclei, introduced in \cite{CEY2013} and shown there to provide an equivalent counting problem with respect to exponential growth.  

Nuclei also constitute a class of triangulations of independent interest. We establish a number of their general properties and provide new examples, including identifying the unique vertex-minimal nontrivial nucleus, with $11$ vertices and $35$ tetrahedra. We also show that every nontrivial nucleus contains a tetrahedron with six interior edges and construct nuclei whose dual graphs have arbitrarily high genus.

Our new reduction (Theorem~\ref{thm:arch-counting}),
together with the barycentric subdivision argument
in the proof of Theorem~\ref{thm:equivalent-exponential-counting},
establishes that simplicial triangulations of $S^3$ have exponential
size if and only if one-vertex generalized
triangulations of $S^3$ do. 
There has been considerable recent interest in these triangulations, both from enumerative and topological perspectives. Altmann and Spreer~\cite{AltmannSpreer2026} used Monte Carlo methods to estimate their number, obtaining numerical evidence that it grows exponentially with the number of tetrahedra, although a rigorous upper bound remains unknown. In another direction, Ibarra, Mathews, Purcell, and Spreer~\cite{IbarraMathewsPurcellSpreer2024} recently proved that every knot in $S^3$ can be represented by a single edge of a one-vertex triangulation of $S^3$.

The paper is organized as follows. In Section~\ref{section:Preliminaries},
we collect relevant definitions and concepts from combinatorics and topology
used throughout the paper. In Section~\ref{section:Enumerate}, we establish
a reduction from simplicial triangulations of $S^3$ to one-vertex generalized
triangulations (Theorem~\ref{thm:arch-counting}) and prove the equivalence
of five exponential counting problems
(Theorem~\ref{thm:equivalent-exponential-counting}).
In Section~\ref{section:Structure}, we investigate
general properties of nuclei and their dual complexes, showing that their
intrinsic polyhedra are closed simple spines of $B^3$ and characterizing when
they are special (Corollary~\ref{cor:special-spine}).
We also introduce a gluing operation that constructs new nuclei from
smaller ones (Theorem~\ref{thm:gluenuclei}).
In Section~\ref{section:GraphGenus}, we study
the dual graphs of nuclei, proving their nonplanarity
(Theorem~\ref{thm:dual-graph-nonplanar}) and constructing
examples whose genus grows linearly with the number of tetrahedra
(Corollary~\ref{cor:linear-genus-growth}).
Finally, Section~\ref{section:Examples} is devoted to examples of nuclei
and their properties, including the identification of the unique vertex-minimal
nontrivial nucleus.

\section{Preliminaries}
\label{section:Preliminaries}
In this section, we collect basic definitions and notation used throughout the paper, in overall agreement with the standard conventions of simplicial and PL topology
\cite{Munkres,RourkeSanderson1982}.

\paragraph{Simplicial complexes and triangulations.}

\begin{figure}[tb]
    \centering
    {\includegraphics[height=9em]{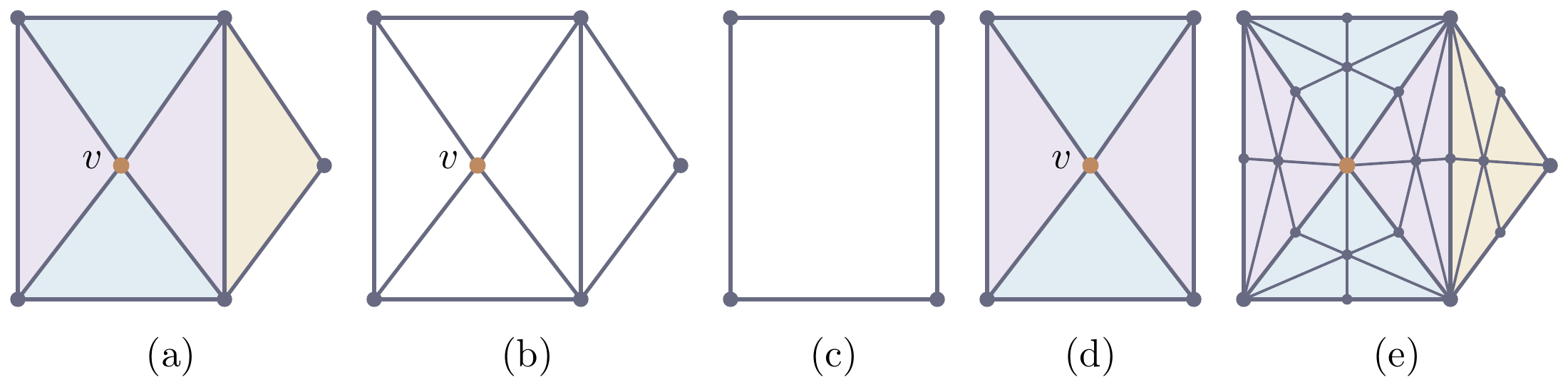}}
    \caption{(a) A simplicial complex $K$; (b) its $1$-skeleton;
    (c) the link of $v$; (d) the star of $v$;
    (e) the barycentric subdivision of $K$.}
    \label{fig:simplicial-definitions}
\end{figure}

\begin{definition}
\label{def:simplicial-triangulation}
A \emph{simplicial complex} $K$ is a finite collection of simplices, 
called its \emph{faces}, such that
\begin{compactenum}
    \item every face of a simplex in $K$ also belongs to $K$;
    \item the intersection of any two simplices in $K$ is a face of both.
\end{compactenum}
The \emph{underlying space} $|K|$ is the union of its simplices.
A \emph{triangulation} of a topological space $X$ is a simplicial
complex $K$ together with a homeomorphism $|K|\to X$.
\end{definition}

Here, the simplices can be thought of as embedded in a sufficiently high-dimensional Euclidean space $\mathbb R^N$, each simplex $\sigma$ being the convex span of its vertices $v_1,\dots, v_{i+1}$, in which case $\sigma$ has \emph{dimension} $i$ and we write $\sigma=[v_1,\dots, v_{i+1}]$. By convention, the empty face has dimension $-1$. We shall not always distinguish between the simplicial complex $K$ as a combinatorial object and its geometric realization in Euclidean space, which is unique up to homeomorphism. 

The inclusion-maximal faces are called \emph{facets}. A simplicial complex is \emph{pure} if all its facets have the same dimension $d$, in which case it is called $d$-dimensional. Its \emph{$k$-skeleton}
consists of all faces of dimension at most $k$; in particular,
the $1$-skeleton consists of the vertices and edges.
Two simplicial complexes are \emph{combinatorially isomorphic} if their face
posets are isomorphic or, equivalently, if they are identical up to a
relabeling of the vertices. We write $f(K)=(f_0,\ldots,f_d)$ for the
\emph{$f$-vector}, where $f_i$ counts the $i$-dimensional faces.

For a triangulated $d$-manifold $K$, the
\emph{boundary} $\partial K$ is the sub-complex generated by the
$(d-1)$-faces belonging to exactly one facet. Its faces are called 
\emph{boundary faces} and $\partial K$ is a triangulation of the boundary of $|K|$; all other faces of $K$ are called \emph{interior}.
For example, if $\Delta^3$ denotes the tetrahedron, then $\partial\Delta^3$ consists of the four triangular
faces of a tetrahedron and all their sub-faces.

\begin{definition}\label{def:links-subdivision}
For a face $\sigma\in K$, the \emph{star} $\operatorname{st}_K(\sigma)$
is the subcomplex consisting of all faces contained in a face of $K$
that contains $\sigma$. Its \emph{simplicial link} is
\[
\operatorname{lk}_K(\sigma)=
\{\tau\in\operatorname{st}_K(\sigma):\tau\cap\sigma=\varnothing\}.
\]
\end{definition}
For example, the link of a vertex in a tetrahedron is the opposite
triangle. It is a useful fact \cite{seifert-threlfall} that a simplicial complex $K$ is a triangulation of a $3$-manifold if every interior vertex has a triangulated $2$-sphere as its link and every boundary vertex has a triangulated $2$-ball as its link.

The \emph{barycentric subdivision} $\operatorname{sd}K$ of a simplicial complex $K$ is defined as usual as having one vertex for each nonempty face of $K$, with simplices corresponding to chains in the poset of faces in $K$. It divides a tetrahedron into $24$ tetrahedra.
These notions are illustrated in Figure~\ref{fig:simplicial-definitions}.

\begin{definition}[Nucleus]\label{def:nucleus}
A \emph{nucleus} is a simplicial triangulation $N$ of $B^3$ such that
all its vertices lie on $\partial N$ and every interior triangle
contains at most one boundary edge. It is \emph{nontrivial} if it
consists of more than one tetrahedron. We denote by $\mathcal N_t$ the set of non-isomorphic nuclei with exactly $t$ tetrahedra.
\end{definition}
The tetrahedron itself is the trivial nucleus. The notion of a nucleus
is closely related to Hachimori's \emph{reduced $3$-balls}
\cite[Definition~4.2]{Hachimori2000}. These satisfy the condition on
interior triangles but may have interior vertices. Thus, nuclei are
precisely reduced $3$-balls without interior vertices.

\begin{figure}[tb]
    \centering
    \includegraphics[height=9em]{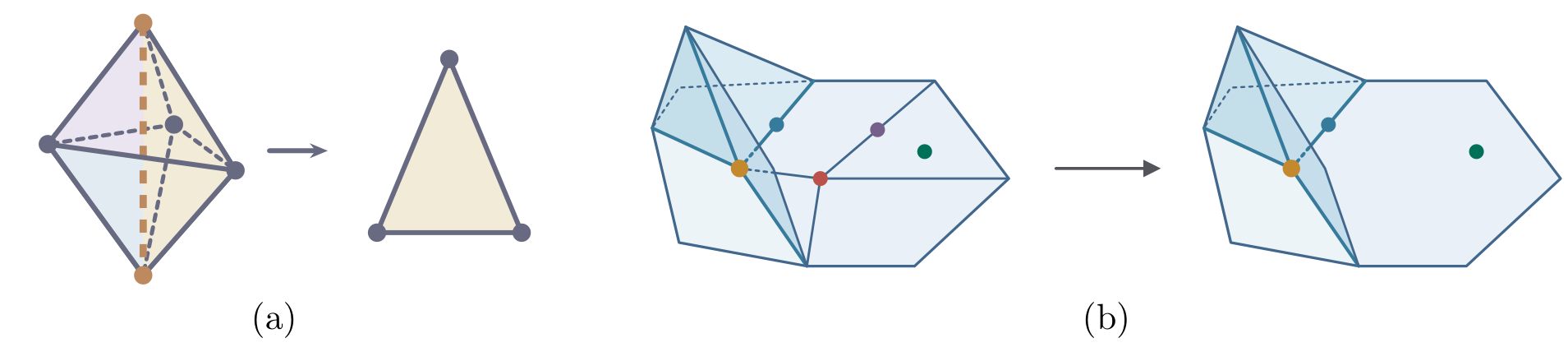}
    \caption{(a) The dual complex of a triangulated 3-ball; (b) A piece of a polyhedron shown with its cellulation and its intrinsic stratification. Each point marked in this figure is a representative of a point with link homeomorphic to one of $S^1$, $\theta$, or $K_4$ }
    \label{fig:dual_Pn}
\end{figure}

\paragraph{Polyhedra, collapses, and dual complexes.}
A \emph{polyhedron} is a space homeomorphic to the geometric
realization of a finite simplicial complex \cite[Ch.~1, \S\S1--2]{Munkres}.
We use the term in this topological sense, rather than
in the convex-geometric or abstract combinatorial sense.
A \emph{CW complex} is a space built by starting with vertices, then attaching edges, followed by disks of increasing dimensions, each along its boundary to cells of lower dimensions. The interiors of these disks are its open cells.
A \emph{cellulation} of a polyhedron $X$ is a finite CW
decomposition of $X$ \cite{Hatcher2002,RourkeSanderson1982}.
Unlike in a simplicial complex, the attaching map of a cell
need not be injective; its boundary may traverse the same
edge more than once. We also use $f$-vectors to record
the numbers of cells in a cellulation.

Following \cite{BL2013}, an \emph{elementary collapse} removes
two faces $\sigma\subsetneq\tau$ of a simplicial complex $K$,
where $\tau$ is the only face that properly contains $\sigma$.
In this case, $\sigma$ is called a \emph{free face};
we write $K\searrow L$ for a sequence of such collapses.
For polyhedra $Y\subseteq X$, we write $X\searrow Y$
if some triangulation of $X$ containing $Y$ as a
subcomplex admits a sequence of elementary
collapses onto $Y$.

\begin{definition}[Dual complex]\label{def:dual-complex}
Let $T$ be a triangulated $3$-ball. Its \emph{dual complex}
$D_T$ is the CW complex with the following cells, whose attaching
maps are induced by the first barycentric subdivision of $T$
\cite{Benedetti2012}:
\begin{compactitem}
    \item Each tetrahedron of $T$ corresponds to a vertex of $D_T$.
    \item Each interior triangle of $T$ corresponds to an edge
    joining the vertices dual to its two tetrahedra.
    \item Each interior edge $e$ of $T$ corresponds to a $2$-cell
    attached along the cycle of dual edges corresponding to the
    triangles containing $e$.
    \item Each interior vertex $v$ of $T$ corresponds to a $3$-cell
    attached along the $2$-sphere formed by the dual $2$-cells
    corresponding to the edges incident to $v$.
\end{compactitem}
\end{definition}

See Figure \ref{fig:dual_Pn}(a) for an illustration. In particular, if $T$ has no interior vertices, as is the case
for nuclei, then $D_T$ has dimension at most $2$. For the trivial
nucleus, $D_T$ consists of one vertex.

For graph terminology, we follow \cite{West2001,BM2001}.
The \emph{dual graph} $G_T$ is the $1$-skeleton of $D_T$:
its vertices correspond to tetrahedra of $T$, with two vertices
adjacent when the corresponding tetrahedra share an interior
triangle. 
The \emph{orientable genus} $\gamma(G)$ of a graph is the
least genus of an orientable surface into which $G$ embeds. We
write $\gamma(N)=\gamma(G_N)$. For example, a planar graph has
genus $0$, while $\gamma(K_5)=1$, where $K_5$ is the complete
graph on five vertices. We also write $C_k$ for the cycle on
$k$ vertices.

\paragraph{Intrinsic polyhedra and spines.}
The \emph{topological link} of a point $x$ in a polyhedron $X$
is the intersection of $X$ with a sufficiently small sphere centered
at $x$ in a geometric realization of $X$. Its homeomorphism type
is independent of these choices.
Let $\theta$ denote the graph with two vertices joined by three edges.

\begin{figure}[tb]
    \centering
    \includegraphics[height=9em]{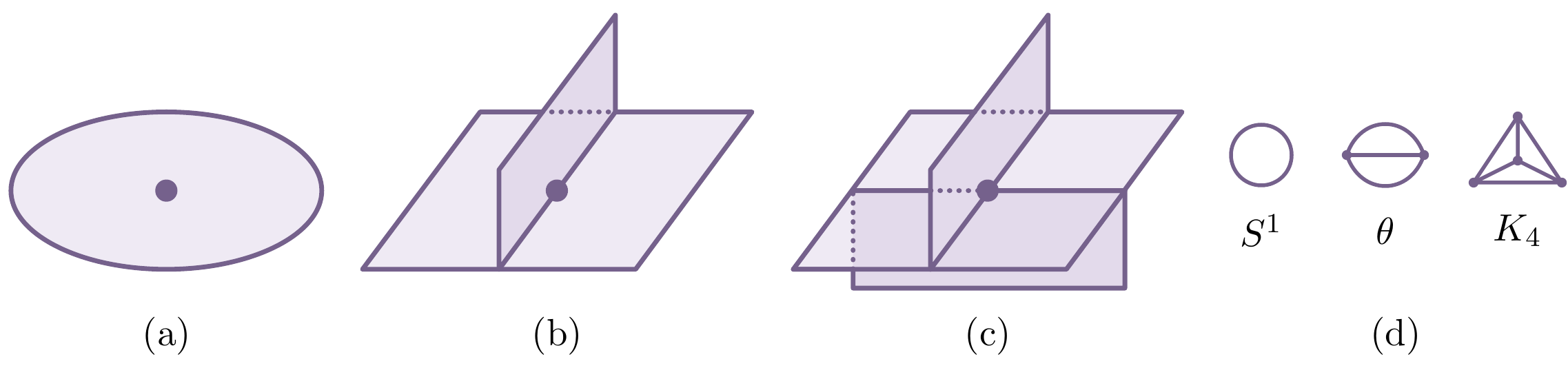}
    \caption{Local models of a closed simple polyhedron:
    (a) an ordinary point; (b) three half-planes meeting along an edge;
    (c) an intrinsic vertex; (d) their respective links
    $S^1$, $\theta$, and $K_4$.}
    \label{fig:simple-polyhedra}
\end{figure}

We use the terminology for closed simple polyhedra and spines from
\cite{Matveev2007,FaganQiuWang2026}.
A compact $2$-dimensional polyhedron is \emph{closed and simple} if every point
has a link homeomorphic to $S^1$, $\theta$, or $K_4$, the complete
graph on four vertices. These are, respectively, the local
models for an ordinary point, a point where three half-planes meet, and an
\emph{intrinsic vertex}. Such polyhedra are also called \emph{closed
fake surfaces}. The local models and their links are shown in
Figure~\ref{fig:simple-polyhedra}.

Following \cite[Definitions~2.2--2.3]{FaganQiuWang2026}, put
\[
\begin{aligned}
S_2(X)&=\{x:\operatorname{lk}_X(x)\cong S^1\},\\
S_1(X)&=\{x:\operatorname{lk}_X(x)\cong\theta\},\\
S_0(X)&=\{x:\operatorname{lk}_X(x)\cong K_4\}.
\end{aligned}
\]
The connected components of $S_2(X)$, $S_1(X)$, and $S_0(X)$
form the \emph{intrinsic stratification} of $X$. They are called
\emph{regions}, \emph{singular edges}, and \emph{intrinsic
vertices}, respectively. The number $c(X)=|S_0(X)|$ is the
\emph{complexity} of $X$.

The intrinsic stratification is \emph{cellular} if its components
are the open cells of a CW decomposition of $X$. In particular, every
region is an open disk and every singular edge is an open interval.
When the stratification is cellular, we simply call \emph{edges of $X$} its singular
edges.
For any finite cell complex $X$, its Euler characteristic is
$\chi(X)=\sum_i(-1)^if_i(X)$; a contractible complex has $\chi(X)=1$. A space is \emph{acyclic} if its reduced homology with integer
coefficients vanishes; in particular, every contractible space is
acyclic.

For a nontrivial nucleus $N$, let $P_N$ denote the underlying
polyhedron $|D_N|$, equipped with its intrinsic stratification.
Figure~\ref{fig:dual_Pn}(b) illustrates the distinction between
an underlying polyhedron's cellulation and intrinsic stratification.

A polyhedron $P\subset M$ is a \emph{spine} of a compact manifold
with boundary if $M\searrow P$. A spine is \emph{special} if it
is closed and simple and its intrinsic stratification is cellular.

\paragraph{Generalized triangulations.}
A \emph{generalized triangulation} of a closed $3$-manifold $M$
is obtained by identifying the triangles of finitely many abstract
tetrahedra in pairs by simplicial maps, so that the resulting space
is homeomorphic to $M$ \cite[Appendix, pp.~534--535]{Hatcher2002}.
Unlike ordinary simplicial triangulations, this allows vertices or
edges of a tetrahedron to be identified. A generalized triangulation
is \emph{one-vertex} if all its vertices are identified. We count
generalized triangulations up to combinatorial equivalence.

\begin{definition}[Exponential class]
A class $\mathcal P$ of combinatorial objects has
\emph{exponential size in terms of a parameter $t$}
if there exists a constant $C\geq 1$ such that,
for all $t$, the number of objects in $\mathcal P$
with the parameter $t$ is at most $C^t$.
Otherwise, $\mathcal P$ has
\emph{superexponential size in terms of $t$}.
\end{definition}

When applying this definition, both the parameter
$t$ and the equivalence relation depend on the
class considered. For instance, we count
triangulations up to combinatorial isomorphism
by their number of tetrahedra, and special spines
up to homeomorphism by their number of
intrinsic vertices.

\section{Counting via dual polyhedra}
\label{section:Enumerate}

We establish the equivalence between exponential counting
problems for nuclei, simplicial triangulations of the
$3$-ball and $3$-sphere, special spines of the $3$-ball,
and one-vertex generalized triangulations of the
$3$-sphere.

\begin{figure}[tb]
    \centering
    \includegraphics[height=10em]{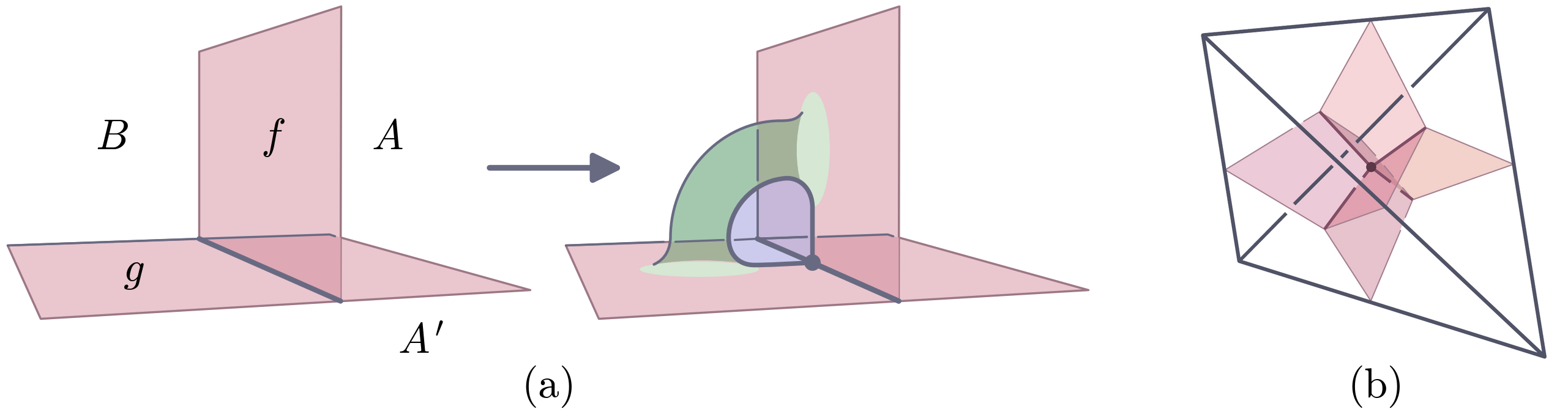}
    \caption{(a) Matveev's arch construction in a polyhedron that connects together the regions $A$ and $A'$; (b) The butterfly of the dual special spine
    inside a tetrahedron, formed by the triangles of
    its first barycentric subdivision that contain no
    original vertex. Redrawn from \cite[Figures~16 and 26]{RST2019}.}
    \label{fig:butterfly}
\end{figure}

The main step is to reduce a simplicial triangulation
of $S^3$ to a one-vertex generalized triangulation,
increasing the number of tetrahedra by at most a factor
of two. We use Matveev's arch construction
\cite[p.~7]{Matveev2007} illustrated in Figure     \ref{fig:butterfly}(a).
Although the construction is naturally expressed using
dual polyhedra, we formulate its consequences directly
for triangulations.

\subsection{From simplicial to one-vertex triangulations}

We begin by recalling the necessary notions and establishing
a key lemma for reducing simplicial triangulations of $S^3$
to one-vertex generalized triangulations.

Let $T$ be a generalized triangulation of a closed
$3$-manifold $M$. We would like to find its spine; however, spines are defined only for manifolds with boundary. Therefore, we first consider the second barycentric subdivision $T^{(2)}$ of $T$, which is a simplicial complex according to \cite[Appendix, p.~535]{Hatcher2002}. We let $M^c$ denote the manifold obtained by removing from $M$ the stars in  $\operatorname{sd}^2 T$ of the original vertices of $T$.  
The union $P$ of all triangles in the first barycentric
subdivision  $\operatorname{sd} T$ that contain no original vertex is a
special spine of $M^c$ \cite[\S2.7]{RST2019}, which we call the
\emph{dual special spine} of $T$.
This distinguishes it from the dual polyhedron $P_N$
of a nucleus, whose intrinsic stratification need not
be cellular (see Remark~\ref{rem:noncellular-nuclei}). Within each tetrahedron of $T$, the dual special
spine forms the butterfly shown in
Figure~\ref{fig:butterfly}(b).

The connected components of $M\setminus P$
correspond bijectively to the vertices of $T$.
Indeed, each tetrahedron of the first barycentric
subdivision has exactly one original vertex, and
its opposite triangle belongs to $P$.
Thus, the component containing an original vertex $v$
is the open cone on its link in $\operatorname{sd}T$.
Since $T$ is a generalized triangulation of a closed
$3$-manifold, every vertex has link homeomorphic to
$S^2$, and hence each component is an open $3$-ball
\cite[\S2.7]{RST2019}.

Each triangle $\tau$ of $T$ corresponds to an singular
edge of $P$, joining the barycenters of its two incident
tetrahedra through the barycenter of $\tau$.
Matveev's arch construction connects two distinct
components of $M\setminus P$ whose closures contain
this edge; see Figure \ref{fig:butterfly}(a).
Moreover, it is proven in \cite[Proposition~5.1]{RST2019} that this
operation preserves the underlying manifold when at
least one of the components of $M\setminus P$ that the arch connects is an open $3$-ball.

\begin{figure}[tb]
    \centering
    \includegraphics[height=10em]
        {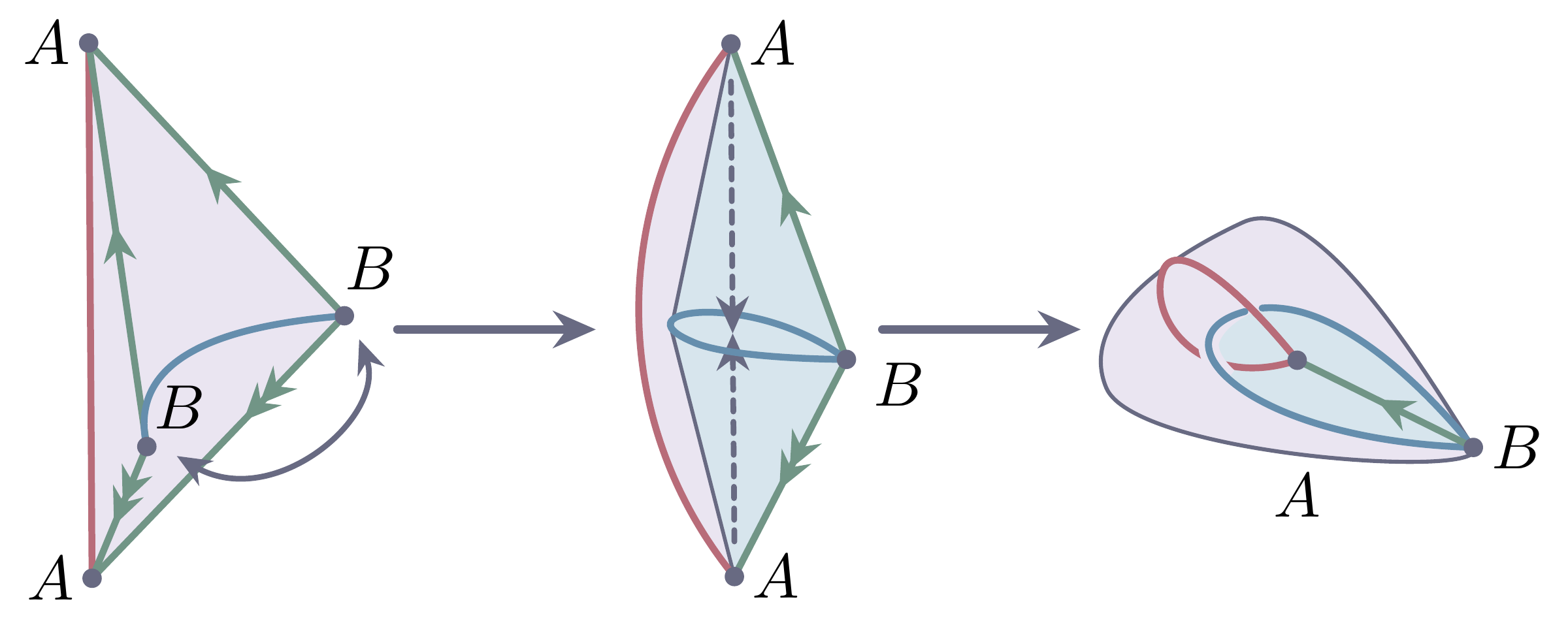}
    \caption{The arch construction in the triangulation.
    First, the pairs of edges marked with matching arrows
    are identified, producing two exterior triangles
    and two cone-shaped interior faces.
    Identifying the interior faces gives a tetrahedron
    with two exterior triangles and two of their vertices
    identified. Redrawn from \cite[Figure~27]{RST2019}.}
    \label{fig:arch}
\end{figure}

In the triangulation, Matveev's arch construction
corresponds to cutting along a triangle $\tau$
and inserting a tetrahedron. Two triangles of the
new tetrahedron are glued to the two exposed copies
of $\tau$, while its remaining two triangles are
identified with each other; see Figure~\ref{fig:arch}.

\begin{lemma}\label{lem:arch-triangulation}
Let $T$ be a generalized triangulation of a closed
$3$-manifold $M$, and let $\tau$ be a triangle containing
two distinct vertices $x\ne y$.
The triangulation analogue of Matveev's arch construction
can be performed at $\tau$ so that the resulting
generalized triangulation $T'$ represents the same
manifold $M$, has one additional tetrahedron, and
has precisely the vertices of $T$ with $x$ and $y$
identified.
\end{lemma}

\begin{proof}
Let $P$ be the dual special spine of $T$.
The vertices $x$ and $y$ correspond to distinct
components of $M\setminus P$, both open $3$-balls
whose closures contain the edge dual to
$\tau$.
Insert an arch connecting these components.
This operation preserves the underlying manifold
\cite[Proposition~5.1]{RST2019}.
Under the correspondence with triangulations,
the operation inserts one tetrahedron into $T$
and identifies precisely $x$ and $y$, leaving
all other vertices distinct; see Figure
\ref{fig:arch}.
This gives the required triangulation $T'$.
\end{proof}

\begin{theorem}\label{thm:arch-counting}
Let $s_t$ and $g_t$ denote the numbers of non-isomorphic
simplicial and one-vertex generalized triangulations of
$S^3$ with $t$ tetrahedra, respectively. Then
\[
s_t\le\sum_{n=t}^{2t-1}g_n.
\]
Consequently, if one-vertex generalized triangulations of
$S^3$ have exponential size, then so do simplicial
triangulations of $S^3$.
\end{theorem}

\begin{proof}
Let $S$ be a simplicial triangulation of $S^3$ with $t$
tetrahedra and $v$ vertices. Since the link of every vertex
is a triangulated $2$-sphere with at least four triangles,
counting vertex-tetrahedron incidences gives $4v\le4t$,
so $v\le t$.

Choose a spanning tree $F$ of the $1$-skeleton of $S$.
We want to assign a distinct triangle $\tau_e$ to every
edge $e\in F$, with $e\subset\tau_e$.
In other words, we seek a system of distinct representatives
for the family $(A_e)_{e\in F}$, where $A_e$ is the set
of triangles containing $e$
\cite[Exercise~3.1.19]{West2001}.

The link of every edge is a cycle of length at least three,
so every edge belongs to at least three triangles. Moreover,
every triangle contains at most two edges of $F$, since $F$
has no cycles. Thus, for every subset $E\subseteq F$,
counting incidences between edges of $E$ and triangles
containing them gives
\[
3|E|\le 2\left|\bigcup_{e\in E}A_e\right|.
\]
In particular,
\[
|E|\le\left|\bigcup_{e\in E}A_e\right|.
\]
This is Hall's condition. By Hall's theorem
\cite[Theorem~3.1.11]{West2001}, the required distinct
triangles $\tau_e$ exist.

We now identify the endpoints of each edge $e\in F$
by applying Lemma~\ref{lem:arch-triangulation} at $\tau_e$.
Since $F$ is a tree, identifying the endpoints of previously
processed edges cannot identify those of an unprocessed
edge $e$, as the unique path between its endpoints in $F$
is $e$ itself.
Moreover, each $\tau_e$ remains available when needed,
since the chosen triangles are distinct and each
tetrahedron insertion changes only the gluing
along its supporting triangle.

After the $v-1$ tetrahedron insertions, all vertices
have been identified. The result is a one-vertex
generalized triangulation $T$ of $S^3$ with
\[
n=t+v-1\le 2t-1
\]
tetrahedra.

To recover $S$ from $T$, we first identify the tetrahedra
introduced by the insertions. These are exactly those
having two of their own triangles identified with each
other. Indeed, no tetrahedron of the original simplicial
triangulation $S$ has such an identification, and
subsequent insertions do not affect previously inserted
tetrahedra.

Since $T$ has $n$ tetrahedra, exactly $v-1=n-t$
of them were introduced by the insertions. Since the insertions were performed at distinct triangles
of $S$, the inserted tetrahedra can be removed independently.
For each removal, the self-identification of the inserted
tetrahedron uniquely determines the gluing of the two
exposed triangles. Thus $S$ is uniquely recovered from $T$.

Consequently,
\[
s_t\le\sum_{n=t}^{2t-1}g_n.
\]
If $g_n\le\alpha^n$ for some $\alpha\ge2$, then
\[
s_t\le\sum_{n=t}^{2t-1}\alpha^n
\le (\alpha^2)^t.
\]
Hence simplicial triangulations of $S^3$ have
exponential size.
\end{proof}

\subsection{Equivalent counting problems}
\label{subsec:equivalent-exponential-counting}
Theorem~\ref{thm:arch-counting} suggests an approach to
the classical exponential-counting problem for triangulated
$3$-spheres: rather than counting triangulations directly,
one may pass to their two-dimensional dual spines and
isolate the combinatorial information that can contribute
more than an exponential factor. Combined with known
correspondences and direct counting arguments, our
reduction establishes the equivalence of the five counting
problems stated in the following theorem.

\begin{theorem}
\label{thm:equivalent-exponential-counting}
The following are equivalent:
\begin{compactenum}
\item the class of nuclei has exponential size;
\item the class of simplicial triangulations of $B^3$ has exponential size;
\item the class of simplicial triangulations of $S^3$ has exponential size;

\item the class of one-vertex generalized triangulations of $S^3$ has exponential size;
\item the class of special spines of $B^3$ has exponential size.
\end{compactenum}
Consequently, these five classes either all have exponential size or all have superexponential size.
\end{theorem}

Let us briefly comment on the individual implications before going into the details of the proof. The equivalence follows from
\[
(1)\Longleftrightarrow(2)\Longleftrightarrow(3)
\Longleftrightarrow(4)\Longleftrightarrow(5).
\]

The nontrivial implication $(1)\Longrightarrow(2)$ is due to
Collet--Eckmann--Younan \cite{CEY2013}, while the reverse implication is immediate because
nuclei form a subclass of simplicial $3$-balls. The equivalence
$(2)\Longleftrightarrow(3)$ is the standard passage between balls and spheres, mentioned for instance in \cite[p.~542]{BB2017}.
The implication $(3)\Longrightarrow(4)$ comes from barycentric subdivision and the
recovery from its dimension-coloring, as in
Rivasseau~\cite[Section~2]{Rivasseau2013SpheresRare}.
The reverse implication $(4)\Longrightarrow(3)$
is our new Theorem~\ref{thm:arch-counting}. The final equivalence
$(4)\Longleftrightarrow(5)$ is the spine--triangulation duality of
Vesnin--Matveev--Fominykh~\cite[Lemma~2.23]{VesninMatveevFominykh2018}. 

\begin{figure}[tb]
    \centering
    \includegraphics[height=9em]{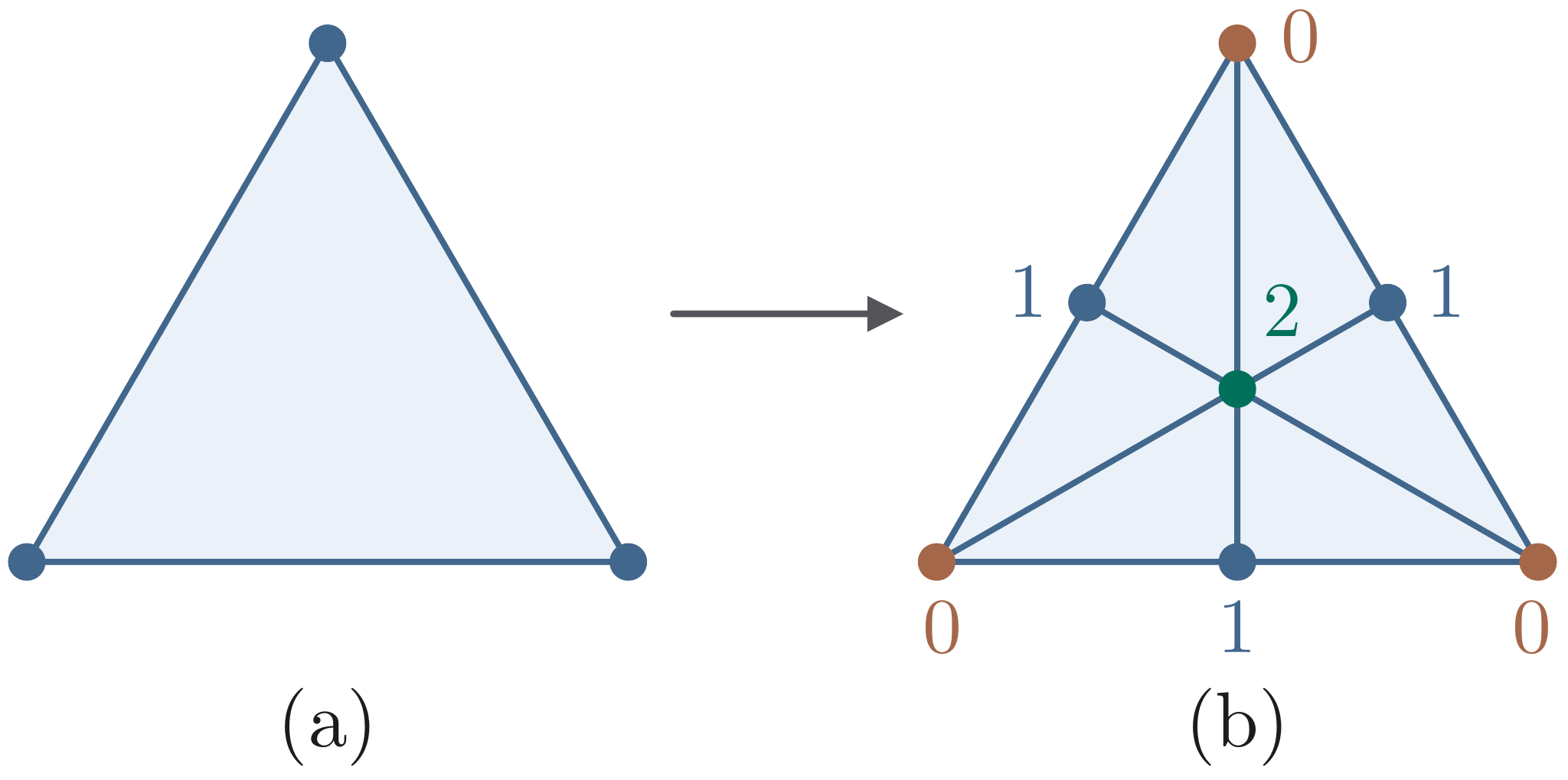}
    \caption{Barycentric subdivision of a triangle.
The vertices are labeled by the dimensions of
the original faces whose barycenters they represent.}
    \label{fig:bar_sub}
\end{figure}

\begin{proof}
\mbox{}\par
\noindent$(1)\Longleftrightarrow(2)$:

\vspace{0.3em}
\noindent For $(1)\Longrightarrow(2)$, Collet--Eckmann--Younan
\cite[Hypothesis~5.6 and Theorem~5.17]{CEY2013}
proved that an exponential bound for rooted nuclei
implies an exponential bound for rooted triangulations
of the $3$-ball.
A nucleus with $t$ tetrahedra has at most $4t$
boundary triangles, each with $3!=6$ possible
orderings of its vertices. Hence it has at most
$24t$ possible roots, so exponential boundedness
of nuclei implies exponential boundedness of
rooted nuclei. 

The
reverse implication $(2)\Longrightarrow(1)$ is immediate because every nucleus is a simplicial
triangulation of $B^3$.
\vspace{1em}

\noindent$(2)\Longleftrightarrow(3)$:

\vspace{0.3em}
\noindent This equivalence is a
standard observation; see, for example,
\cite[p.~542]{BB2017}. For completeness, we
include a short proof.

For $(2)\Longrightarrow(3)$, remove a tetrahedron
from a simplicial triangulation of $S^3$.
This gives a simplicial triangulation of $B^3$
with one fewer tetrahedron, whose boundary
consists of the four triangles of the removed
tetrahedron. Since the original sphere is uniquely
recovered by gluing the tetrahedron back, an
exponential bound for simplicial $3$-balls gives
one for simplicial $3$-spheres.

For $(3)\Longrightarrow(2)$, let $B$ be a simplicial
triangulation of $B^3$ with $t$ tetrahedra.
For a new vertex $w$, the cone $w*\partial B$
consists of $\partial B$ together with all simplices
obtained by adjoining $w$ to its faces.
Gluing this cone to $B$ along $\partial B$ produces
a simplicial triangulation of $S^3$.
Since the boundary of $B$ has at most $4t$ triangles,
the resulting sphere has at most $5t$ tetrahedra.

The original ball is recovered by deleting the cone
vertex. A simplicial $3$-sphere with $m$ tetrahedra
has at most $m$ vertices, as shown in the proof of
Theorem~\ref{thm:arch-counting}.
Hence, if $s_m\le K^m$ for some $K\ge1$, then
the number of simplicial $3$-balls with $t$
tetrahedra is at most
\[
\sum_{m\le5t}m\,s_m
\le25t^2K^{5t}.
\]

\vspace{1em}

\noindent $(3)\Longleftrightarrow(4)$:

\vspace{0.3em}
\noindent The implication $(4)\Longrightarrow(3)$ is
Theorem~\ref{thm:arch-counting}.

For $(3)\Longrightarrow(4)$, let $T$ be a one-vertex generalized triangulation
of $S^3$ with $t$ tetrahedra.
As mentioned above~\cite[Appendix, p.~535]{Hatcher2002},
two barycentric subdivisions of a
generalized triangulation produce a simplicial
triangulation. In dimension $3$, each subdivision
replaces every tetrahedron by $24$ tetrahedra.
Thus $\operatorname{sd}^2 T$ is a simplicial
triangulation of $S^3$ with $576t$ tetrahedra.

Rivasseau~\cite[Section~2]{Rivasseau2013SpheresRare}
explained how a triangulation can be recovered
from its barycentric subdivision using the
coloring of vertices by the dimensions of the
original faces (see Figure \ref{fig:bar_sub}). The same reconstruction applies
to generalized triangulations. Indeed, each vertex of color $3$ corresponds to an
original tetrahedron, whose barycentric subdivision
consists of the $24$ tetrahedra incident to this vertex.
Each of its four triangles is subdivided into six
triangles, and their adjacencies determine the original
triangle identifications, including those between triangles
of the same tetrahedron.
Every tetrahedron of a barycentric subdivision
contains exactly one vertex of each of the four
colors. Since the dual graph is connected, the
coloring is determined by its values on one
tetrahedron. Hence there are at most $4!=24$
possible colorings, each determining at most
one original triangulation.

Applying this reconstruction to both subdivisions,
a fixed simplicial triangulation has at most
$24^2=576$ generalized triangulations as preimages.
Therefore, if $s_t$ denotes the number of simplicial
triangulations of $S^3$ with $t$ tetrahedra and
$g_t$ denotes the number of one-vertex generalized
triangulations of $S^3$ with $t$ tetrahedra, then
\[
g_t\le 576\,s_{576t}.
\]
Consequently, if $s_t\le K^t$ for some constant
$K>0$, then $g_t$ is exponentially bounded in $t$:
\[
g_t\le 576\,K^{576t}.
\]

\vspace{1em}

\noindent$(4)\Longleftrightarrow(5)$:

\vspace{0.3em}
\noindent  
Vesnin--Matveev--Fominykh
\cite[Lemma~2.23]{VesninMatveevFominykh2018}
proved that special spines of a closed $3$-manifold $M$
correspond bijectively to its one-vertex generalized
triangulations, with intrinsic vertices corresponding
to tetrahedra. Here a special spine of a closed manifold $M$ means a special spine of $M$ with an open $3$-ball removed,
as in the dual construction discussed above
for one-vertex triangulations. Let $b_t$ denote the number of homeomorphism classes
of special spines of $B^3$ with $t$ intrinsic vertices.
Since removing an open $3$-ball from $S^3$ leaves
$B^3$, in our setting this correspondence gives
$g_t=b_t$.
\end{proof}

At present, an exponential upper bound is not known for any of the five
equivalent classes. The strongest direct evidence comes from one-vertex
generalized triangulations of $S^3$. Using Monte Carlo sampling,
Altmann--Spreer experimentally observe that their number appears to grow
singly exponentially with the number of tetrahedra
\cite{AltmannSpreer2026}. This is computational evidence rather than a rigorous
upper bound.

On the dual side, Fagan--Qiu--Wang obtain a complete classification of
contractible cellular simple polyhedra through complexity four. At
complexities $1,2,3,4$ they find respectively
\[
2,\ 17,\ 238,\ 4618
\]
such polyhedra, of which
\[
1,\ 3,\ 20,\ 128
\]
are special spines~\cite[Table~2.1]{FaganQiuWang2026}. These data suggest that
special spines become increasingly sparse, and they conjecture that their
proportion tends to zero as the complexity tends to infinity
\cite[Conjecture~2.5]{FaganQiuWang2026}. This conjecture concerns the relative
abundance of special spines inside the larger class of contractible cellular
simple polyhedra, and does not by itself predict whether the number of special
spines is exponential or superexponential.

\section{Structural properties of nuclei}
\label{section:Structure}
The counting equivalences established above motivate
a closer examination of nuclei and their dual complexes.
The condition defining a nucleus is local: it restricts how many
boundary edges an interior triangle may contain. Nevertheless,
after passing to the dual complex this local condition has
rather strong global consequences. In this section, we show
that the dual complex of every nontrivial nucleus is a
contractible closed simple polyhedron and that its intrinsic
polyhedron is a closed simple spine of $B^3$. We also characterize
when this spine is special and introduce a gluing operation
that will be useful when studying the genus of dual graphs
in Section~\ref{section:GraphGenus}.

\subsection{The dual complex of a nucleus}
\label{subsec:dual-complex}
We now study the dual complex $D_N$ of a nontrivial nucleus $N$.
Recall from Section~\ref{section:Preliminaries} that, since all vertices of a
nucleus lie on the boundary, $D_N$ has no $3$-cells and is therefore
$2$-dimensional. Its $0$-cells correspond to tetrahedra of $N$, its
$1$-cells to interior triangles, and its $2$-cells to interior edges.
We first record a simple consequence of the nucleus condition.

\begin{lemma}
\label{lem:one-boundary-triangle}
Let $\sigma=abcd$ be a tetrahedron of a nontrivial nucleus $N$. Then
$\sigma$ contains at most one boundary triangle. Moreover, if $abc$ is a
boundary triangle of $\sigma$, then the three edges
\(
ad,\  bd,\  cd
\)
are interior.
\end{lemma}

\begin{proof}
Let $\sigma=abcd$ be a tetrahedron of a nontrivial nucleus $N$. Suppose first that $\sigma$ contains two boundary triangles, say $abc$ and
$abd$. Then
\[
ab,ac,bc,ad,bd
\]
are boundary edges. If $acd$ were an interior triangle, it would contain the
two boundary edges $ac$ and $ad$, which contradicts the nucleus condition.
Hence $acd$ is a boundary triangle. Similarly, $bcd$ is a boundary triangle. Since $N$ is a triangulation of a ball, $N=\sigma$, contrary to the assumption that $N$ is nontrivial. Therefore, $\sigma$ has at most one boundary triangle.

Now suppose that $abc$ is a boundary triangle. By the first part,
$abd,acd,bcd$ are interior triangles. The triangle $abd$ already contains
the boundary edge $ab$, so the nucleus condition forces $ad$ and $bd$ to be
interior. Similarly, $cd$ is interior.
\end{proof}

The first structural result shows that the nucleus condition forces $D_N$ to have precisely the local topology of a closed simple polyhedron.

\begin{theorem}
\label{thm:dual-closed-simple}
Let $N$ be a nontrivial nucleus. Then $D_N$ is a compact, contractible,
closed simple $2$-dimensional polyhedron.
\end{theorem}

\begin{proof}
Since $N$ is finite, $D_N$ is compact. Moreover,
$N$ is a triangulation of the $3$-ball, hence a
PL $3$-ball \cite[Theorem~2.3]{Schleimer2007}.
The dual complex of a PL ball is homotopy
equivalent to the ball
\cite[Lemma~3.2]{Benedetti2012}.
Thus $D_N$ is contractible.

It remains to understand the local structure of $D_N$.
\begin{compactitem}
    \item A point in the interior of a dual $2$-cell has a disk neighborhood, and hence link $S^1$.
    \item Let $x$ lie instead in the interior of the dual edge $\tau^*$, where $\tau$ is
    an interior triangle of $N$. The dual $2$-cells incident with $\tau^*$
    correspond exactly to the interior edges of $\tau$. Since $N$ is a nucleus,
    $\tau$ has either two or three interior edges. Accordingly, two or three 2-cells meet along
    $\tau^*$. Thus
    \[
    \operatorname{lk}_{D_N}(x)\cong S^1
    \qquad\text{or}\qquad
    \operatorname{lk}_{D_N}(x)\cong\theta.
    \]
    In particular, $D_N$ has no boundary.
    \item It remains to consider the dual vertices. Let
    \(
    \sigma=abcd
    \)
    be a tetrahedron of $N$, and let $v_\sigma$ be the corresponding vertex of
    $D_N$. For every interior triangle $\tau\subset\sigma$, let $x_\tau$ denote the
    corresponding vertex of
    \[
    L_\sigma:=\operatorname{lk}_{D_N}(v_\sigma).
    \]
    If two interior triangles $\tau,\tau'\subset\sigma$ meet in an interior edge
    $e$, then the dual $2$-cell $e^*$ contributes the edge
    $x_\tau x_{\tau'}$ to $L_\sigma$. Hence
    \begin{equation}
    \label{eq:dual-link-dictionary}
    \begin{aligned}
    V(L_\sigma)
      &=\{x_\tau:\tau\subset\sigma\text{ is an interior triangle}\},\\
    E(L_\sigma)
      &=\{x_\tau x_{\tau'}:
          \tau\cap\tau'\text{ is an interior edge of }N\}.
    \end{aligned}
    \end{equation}
    
    By Lemma~\ref{lem:one-boundary-triangle}, $\sigma$ contains either one boundary triangle or none.    
    \begin{compactitem}
        \item Suppose first that $\sigma$ has one boundary triangle, say $abc$. Then
        $abd,acd,bcd$ are interior triangles and $ad,bd,cd$ are interior edges.
        By \eqref{eq:dual-link-dictionary},
        \[
        L_\sigma\cong C_3\cong S^1.
        \]
         \item Suppose next that $\sigma$ has no boundary triangle. Two boundary edges of $\sigma$ cannot share a
        vertex, since the interior triangle containing them would contain two boundary
        edges. There are
        therefore only three possibilities: 0, 1, or 2 boundary edges, and in the last case, the two boundary edges are opposite.        
        Using \eqref{eq:dual-link-dictionary}, the corresponding links of $v_\sigma$ are
        respectively
        \[
        K_4,\qquad K_4-e\cong\theta,\qquad C_4\cong S^1.
        \]
    \end{compactitem}     
\end{compactitem}
Therefore, $D_N$ is a compact, contractible, closed simple $2$-dimensional polyhedron.
\end{proof}

    The proof gives a complete description of the possible links of dual
    vertices:
    \[
    \begin{array}{c|c|c}
    \text{boundary triangles of }\sigma
    &
    \text{boundary edges of }\sigma
    &
    \operatorname{lk}_{D_N}(\sigma^*)\ 
    \\ \hline
    1 & 3 & C_3\ \\
    0 & 2\text{ opposite} & C_4\ \\
    0 & 1 & K_4-e\ \\
    0 & 0 & K_4\ 
    \end{array}
    \]
    In particular,
    \begin{equation}
    \label{eq:k4-all-interior}
    \operatorname{lk}_{D_N}(v_\sigma)\cong K_4
    \quad\Longleftrightarrow\quad
    \text{all six edges of }\sigma\text{ are interior}.
    \end{equation}
    This last observation will be the key to the next subsection.

\subsection{An all-interior-edge tetrahedron}
\label{subsec:all-interior-tetrahedron}

The local analysis above shows that a $K_4$-link is equivalent to a
tetrahedron with all six edges interior. We now show that this local
configuration must occur in every nontrivial nucleus.

\begin{theorem}
\label{thm:all-interior-tetrahedron}
Every nontrivial nucleus contains a tetrahedron all six of whose edges are
interior.
\end{theorem}

\begin{proof}
By Theorem~\ref{thm:dual-closed-simple}, \(D_N\) is a contractible closed simple polyhedron, and hence it is acyclic.

A result of Ikeda~\cite{Ikeda1971}, in the formulation of Koda and Naoe~\cite[Lemma~2.7(2), p.~3718]{KodaNaoe2020}, implies that an acyclic closed simple polyhedron has a point whose link is isomorphic to \(K_4\). Since \(D_N\) is closed, it follows that \(D_N\) contains a point with \(K_4\)-link.

By local analysis in the proof of Theorem~\ref{thm:dual-closed-simple}, such a point is a dual vertex \(\sigma^*\). By \eqref{eq:k4-all-interior}, the corresponding tetrahedron \(\sigma\) has all six edges interior.
\end{proof}

This theorem provides an elegant obstruction to planarity
of the dual graph of a nontrivial nucleus, as shown in
Section~\ref{section:GraphGenus}.

\subsection{The intrinsic polyhedron \texorpdfstring{$P_N$}{PN}}
\label{subsection:pn} 
Recall that $P_N$ is the underlying polyhedron of $D_N$,
with its intrinsic stratification. Thus, a dual vertex
with link $C_3$ or $C_4$ is an ordinary point of a region,
one with link $K_4-e$ lies on a singular edge, and one
with link $K_4$ is an intrinsic vertex.

We introduce three notions that will be used below.
Let $X$ be a contractible closed simple polyhedron.
Its \emph{singular set} $S=S_0(X)\cup S_1(X)$ consists
of the points where $X$ is not locally a surface.
Let $S^{(1)},\ldots,S^{(k)}$ be the connected components
of $S$. Let $U_1,\ldots,U_k$ be the interiors of pairwise
disjoint regular neighbourhoods of these components,
as in \cite[Ch.~3]{RourkeSanderson1982}, chosen small
enough that removing their union does not disconnect
any region. Put $U=\bigcup_i U_i$.

For each region $R$, the \emph{truncated region}
$F_R=R\setminus U$ is a compact connected surface
with boundary. For example, if $R$ is an open disk,
then $F_R$ is a closed disk with one boundary circle.
If $R$ is an open annulus, then $F_R$ is a closed
annulus with two boundary circles. Let $h_R$ denote
the number of boundary circles of $F_R$.

The \emph{incidence graph} of $X$ has one vertex for
each region $R$ and each component $S^{(i)}$ of $S$.
Each boundary circle of $F_R$ lies on the boundary
of exactly one $U_i$ and determines an edge joining
$R$ to $S^{(i)}$. Thus the graph is bipartite,
and multiple edges are allowed.

\begin{lemma}\label{lem:cellular-criterion}
Let $X$ be a contractible closed simple polyhedron
with $c$ intrinsic vertices, $r$ regions, and
$k$ components of its singular set $S$.
Then the following hold.

\begin{compactenum}
\item The number of regions is $r=c+k$,
every truncated region $F_R$ is planar,
and the incidence graph is a tree.

\item The intrinsic stratification of $X$ is cellular
if and only if $S$ is connected and contains an
intrinsic vertex. In this case, $X$ has $2c$
edges and $c+1$ regions.
\end{compactenum}
\end{lemma}

\begin{proof}
\emph{For 1.}, by a classical result of Ikeda~\cite{Ikeda1971},
see also \cite[Lemma~2.7(1)]{KodaNaoe2020},
every region of an acyclic closed simple polyhedron is
a sphere with holes. Since $X$ is contractible,
it is acyclic, so Ikeda's result applies.
Consequently, every truncated region $F_R$ is
a sphere with holes and hence planar.

We next show that the incidence graph $\Gamma$
is a tree. By \cite[Lemma~2.8]{KodaNaoe2020},
if $C$ is a simple closed curve contained in a
region of an acyclic closed simple polyhedron, then
removing $C$ disconnects the polyhedron.
This applies to every boundary circle $C$ of a
truncated region $F_R$. If the corresponding edge
of $\Gamma$ belonged to a cycle, the remaining
edges of that cycle would provide another path
between the two sides of $C$, contradicting the
separation property. Hence removing any edge
disconnects $\Gamma$. The graph $\Gamma$ is
connected because $X$ is connected. Therefore, $\Gamma$ is a tree.

We now count the regions using the Euler characteristic
of finite CW complexes. As for simplicial complexes,
it is the alternating sum of the numbers of cells,
and it is invariant under homotopy equivalence
\cite[Theorem~2.44]{Hatcher2002}.

The singular set $S$ is a finite graph whose intrinsic
vertices have degree four. It may also have components
homeomorphic to circles, containing no intrinsic vertices.
Let $a$ be the number of such circular components.
By introducing one auxiliary vertex on each circle,
we obtain a graph with $c+a$ vertices and
$2c+a$ edges. Hence
\[
\chi(S)=(c+a)-(2c+a)=-c.
\]

Let $B=\sum_R h_R$ be the total number of boundary
circles of the truncated regions. The incidence
graph has $r+k$ vertices and $B$ edges. Since
it is a tree, $B=r+k-1$.

Each $F_R$ is planar, so $\chi(F_R)=2-h_R$.
Since the closure $\overline{U}$ of $U$ retracts
onto $S$, we have $\chi(\overline{U})=\chi(S)$. Moreover,
$\overline{U}\cap F_R=\partial F_R$ is a union of $h_R$
circles, each with Euler characteristic zero.
Applying the additivity of Euler characteristic
\cite[\S 2.2, Exercise~21]{Hatcher2002}, we obtain $r=c+k$ by
\[
1=\chi(X)
=\chi(\overline{U})+\sum_R
  \bigl(\chi(F_R)-\chi(\partial F_R)\bigr)
=\chi(S)+\sum_R\chi(F_R)
=-c+2r-B
=-c+r-k+1.
\]

\emph{For 2.}, by \cite[Lemma~2.7(2)]{KodaNaoe2020},
every  acyclic closed simple polyhedron has at
least one intrinsic vertex.
Suppose first that the intrinsic stratification
of $X$ is cellular. Then every truncated region
$F_R$ is a disk and hence has exactly one boundary
circle. Thus each region is adjacent to only one
component of $S$ in the incidence graph. If $S$
had more than one component, the incidence graph
would therefore be disconnected, a contradiction.

Conversely, suppose that $S$ is connected and contains an intrinsic vertex.
The incidence graph then has exactly one vertex
corresponding to the single component of $S$.
Since it is a tree, every region is joined to
this vertex by exactly one edge, so $h_R=1$
for every region $R$. Since each $F_R$ is also
planar, it is a disk. Moreover, $S$ is connected
and contains an intrinsic vertex, so every
singular edge is an open interval. Hence the
intrinsic stratification is cellular.

Finally, since $S$ is connected, it is a
$4$-regular graph with $c$ vertices. By the
degree-sum formula,
the number $e$ of its edges satisfies
$2e=4c$, so $e=2c$.
Moreover, $k=1$, and part~(1) gives
$r=c+k=c+1$.
\end{proof}

The singular set of the dual polyhedron $P_N$ can be
described directly in terms of the dual graph $G_N$,
using the link table following
Theorem~\ref{thm:dual-closed-simple}.
The following proposition gives this description
and characterizes when the intrinsic stratification
of $P_N$ is cellular.

\begin{proposition}\label{prop:intrinsic-structure}
Let $N$ be a nontrivial nucleus, and let $c$ be the
number of intrinsic vertices of $P_N$.
Then the following hold.

\begin{compactenum}
\item The singular set of $P_N$ consists of the
edges of $G_N$ dual to interior triangles whose
three edges are interior, together with their
endpoints. Its intrinsic vertices correspond
exactly to tetrahedra with all six edges interior.
In particular, $c\ge1$.

\item The intrinsic stratification of $P_N$ is
cellular if and only if its singular set is
connected. In this case, it has $2c$ edges
and $c+1$ regions, all open disks.
\end{compactenum}
\end{proposition}

\begin{proof}
\emph{For 1.}, by Theorem~\ref{thm:all-interior-tetrahedron}
and \eqref{eq:k4-all-interior}, every nontrivial
nucleus contains a tetrahedron with all six edges
interior. By the link table following
Theorem~\ref{thm:dual-closed-simple}, these
tetrahedra correspond exactly to the intrinsic
vertices of $P_N$. Hence $c\ge1$.
The same link table shows that a dual edge is
singular precisely when its corresponding interior
triangle has three interior edges. Its endpoints
account for all remaining singular points.
\emph{For 2.}, by Theorem~\ref{thm:dual-closed-simple},
$P_N$ is a contractible closed simple polyhedron.
The claim now follows directly from
Lemma~\ref{lem:cellular-criterion}.
\end{proof}

The dual polyhedron is a spine even when its
intrinsic stratification is not cellular.

\begin{corollary}\label{cor:special-spine}
For every nontrivial nucleus $N$, $P_N$ is a closed
simple spine of $B^3$. It is a special spine
if and only if its singular set is connected.
In that case, we call $N$ \emph{special}.
\end{corollary}

\begin{proof}
Benedetti showed that every triangulated PL $3$-ball
$T$ collapses onto its dual complex $D_T$
\cite[Lemma~3.2]{Benedetti2012}.
Thus $|N|\searrow |D_N|$, so $P_N$ is a spine
of $B^3$. By Proposition~\ref{prop:intrinsic-structure},
$P_N$ is special if and only if its singular set
is connected.
\end{proof}

A dual polyhedron can therefore fail to be special even though it is a
spine of the $3$-ball. The following gluing construction gives concrete
examples of this phenomenon, and will also be useful in our study of dual
graphs.

\subsection{Gluing nuclei}
\label{subsec:gluing-nuclei}

\begin{figure}[tb]
    \centering
    \includegraphics[height=9em]{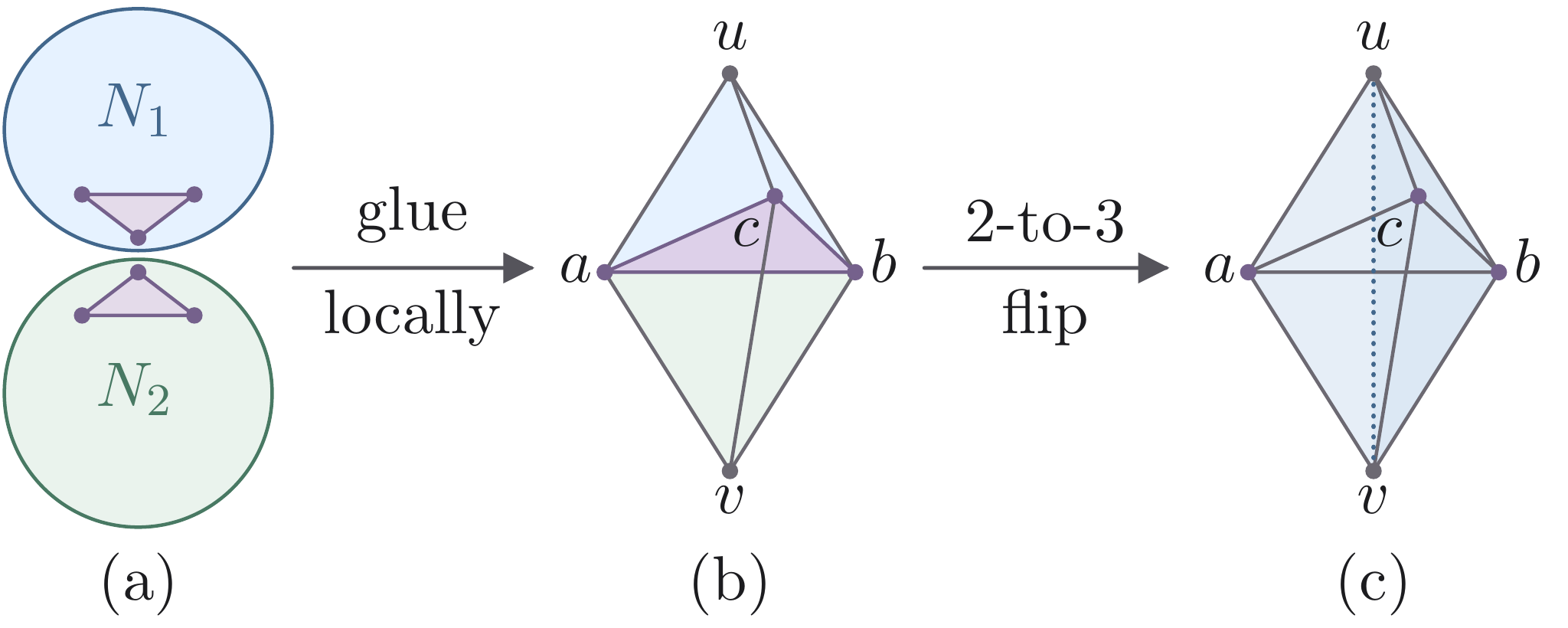}
    \caption{Gluing two nuclei:
    (a) the chosen boundary triangles;
    (b) the tetrahedra after identifying these triangles;
    (c) the $2$-to-$3$ bistellar move introducing
    the interior edge $uv$.}
    \label{fig:gluing-nuclei}
\end{figure}

We finish this section with an operation that combines two nuclei without leaving the class. Its definition is local, but later it will allow us to build nuclei
with increasingly complicated dual graphs.

Given two disjoint nuclei, $N_1$ and $N_2$, and
boundary triangles $\Delta_1$ in $N_1$ and
$\Delta_2$ in $N_2$, we denote by
\[
{N_1}_{\Delta_1}\sharp_{\Delta_2}N_2
\]
the triangulation obtained by identifying
$\Delta_1$ and $\Delta_2$ and performing a
$2$-to-$3$ bistellar flip on the resulting
interior triangle (see Figure \ref{fig:gluing-nuclei}). More precisely, if $t_1$
and $t_2$ are the tetrahedra containing
$\Delta_1$ and $\Delta_2$, respectively,
we replace $t_1\cup t_2$ by three tetrahedra
sharing the new edge $(x_1\,x_2)$, where
$x_i$ is the vertex of $t_i$ opposite to
$\Delta_i$, for $i=1,2$.

\begin{theorem}
\label{thm:gluenuclei} 
If at least one of the nuclei $N_1$ and $N_2$ is not a tetrahedron, then ${N_1}_{\Delta_1}\sharp_{\Delta_2} N_2$ is a nucleus.
\end{theorem}

\begin{proof} It is clear by construction that ${N_1}_{\Delta_1}\sharp_{\Delta_2} N_2$ is a triangulation of the ball $B^3$ with all vertices on the boundary.
Moreover, the interior triangles of ${N_1}_{\Delta_1}\sharp_{\Delta_2} N_2$ coincide with those of $N_1\cup N_2$ except that $\Delta_1=\Delta_2$ in the latter is replaced by the three interior triangles of $t_1'\cup t_2'\cup t_3'$, which are of the form $\Delta = (x_1\, x_2\, y)$, where  $y$ belongs to $\Delta_1 =\Delta_2$, and obviously $(x_1\, x_2)$ is an interior edge in ${N_1}_{\Delta_1}\sharp_{\Delta_2} N_2$. On the other hand,  the boundary of ${N_1}_{\Delta_1}\sharp_{\Delta_2} N_2$ is identical to the boundary of $N_1\cup N_2$, so if $(x_1\, y)$ were a boundary edge for some vertex $y$ in $\Delta_1 =\Delta_2$, then this would be the case for all $y$ in $\Delta_1 =\Delta_2$, since $N_1$ is a nucleus. This would imply that $N_1$ equals a single tetrahedron. Similarly, if $(x_2\, y)$ is a boundary edge for some vertex $y$ in $\Delta_1 =\Delta_2$, then $N_2$ is a single tetrahedron. Thus, if either $N_1$ or $N_2$ is not a tetrahedron, then at most one of $(x_1\,y)$ and $(x_2\, y)$ is a boundary edge for each vertex $y$ in $\Delta_1 =\Delta_2$. This completes the proof of the theorem. 
\end{proof}  

When the gluing triangles are not relevant, we write $N_1\sharp N_2$
instead of ${N_1}_{\Delta_1}\sharp_{\Delta_2}N_2$. Clearly, the hypothesis of the theorem is necessary: if both factors are single tetrahedra, then
each new interior triangle $x_1x_2y$ contains the two boundary edges
$x_1y$ and $x_2y$.

We note, as a first consequence of this gluing property, the following result on the asymptotic behaviour of the number $|{\mathcal N}_t|$ of nuclei of size $t$.

\begin{proposition} The limit 
$$\kappa := \lim_{t\to\infty} \frac{\ln|{\mathcal N}_t|}{t}$$
exists in $\mathbb R_+\cup\{\infty\}$. If it is finite, it holds for large enough $t$ that
$$ |{\mathcal N}_t|\,\leq\,  e^{\kappa(t+1)}\,.$$    
\end{proposition}

\begin{proof} Let $\tilde{\mathcal N}_t$ denote the set of (non-isomorphic) nuclei of size $t$ and with one marked oriented edge on the boundary. Given $T\in\tilde{\mathcal N}_t$ and $S\in\tilde{\mathcal N}_s$, where $t,s\geq 2$, with marked oriented edges $e_T$ and $e_S$, respectively, let $\Delta_1$ be the boundary triangle of $T$ to the left of $e_T$ (relative to the standard orientation of ${\partial B}^3$) and let $\Delta_2$ be the boundary triangle of $S$ to the right of $e_S$. We then denote by $T\sharp S$ the nucleus in $\tilde{\mathcal N}_{t+s+1}$ obtained by gluing $T$ to $S$ along $\Delta_1$ and $\Delta_2$, such that $e_T$ and $e_S$ get identified and the resulting oriented edge in its boundary is chosen as the marked edge. Since $T$ and $S$ can be uniquely reconstructed from $T\sharp S$, it follows that
$$|\tilde{\mathcal N}_{t+s+1}| \geq |\tilde{\mathcal N}_{t}|\,|\tilde{\mathcal N}_{s}|\,.$$
Because there are nuclei with $q$ tetrahedra for every $q\ge 37$ (see Remark \ref{rem:nuclei-all-sizes}), the above shows that the function $t\to \ln|\tilde{\mathcal N}_{t-1}|$ is super-additive for large enough $t$:
$$\ln|\tilde{\mathcal N}_{t+s-1}| \geq \ln|\tilde{\mathcal N}_{t-1}|\,+\,\ln|\tilde{\mathcal N}_{s-1}|\,.$$
It is a standard consequence, using Fekete's lemma, that
$$  \lim_{t\to\infty}\frac{\ln|\tilde{\mathcal N}_{t-1}|}{t} = \sup_{t\geq 38}\frac{\ln|\tilde{\mathcal N}_{t-1}|}{t}\,.$$
 As there are at most $t$ triangles in the boundary of $T\in{\mathcal N}_t$, there are at most $3t$ ways of orienting an edge in the boundary, which implies that
 $3t|{\mathcal N}_t|\geq |\tilde{\mathcal N}_t|$. Together with the obvious inequality $|{\mathcal N}_t|\leq |\tilde{\mathcal N}_t|$ and the preceding identity, we conclude that the limit $\kappa$ exists and that 
 $$\kappa = \lim_{t\to\infty} \frac{\ln|{\mathcal N}_t|}{t} = \lim_{t\to\infty} \frac{\ln|\tilde{\mathcal N}_t|}{t} = \sup_{t\geq 38}\frac{\ln|\tilde{\mathcal N}_{t-1}|}{t}\geq \sup_{t\geq 38}\frac{\ln|{\mathcal N}_{t-1}|}{t}\,,$$
which also proves the last statement of the proposition, if $\kappa$ is finite. 
\end{proof}

We will return to gluing in Section~\ref{section:GraphGenus}, where it
produces nuclei of unbounded dual genus. The following remark
shows that the intrinsic polyhedra of nuclei need not be special.

\begin{remark}\label{rem:noncellular-nuclei}
The cellularity assertion of Fagan, Qiu, and Wang
\cite[Proposition~2.4(1)]{FaganQiuWang2026}
needs an additional hypothesis. For a simple counterexample,
take two copies of the \emph{abalone}, a special spine of $B^3$
with one intrinsic vertex \cite[Table~2.1]{FaganQiuWang2026},
and identify a small closed disk in the interior of a region
of each copy. The resulting polyhedron is contractible,
since the two copies and their intersection are contractible
\cite[Chapter~0, Exercise~23]{Hatcher2002}.
It is still closed and simple: along the boundary of
the identified disk, three half-planes meet. This boundary is an isolated
singular circle with no intrinsic vertices, and the two punctured
regions are annuli. Hence its intrinsic stratification is not
cellular.

The same obstruction occurs for nuclei. Let $N_1,N_2$ be nontrivial
nuclei, let $\Delta=[y_1,y_2,y_3]$ be the identified boundary
triangle, and let $x_i$ be the vertex opposite $\Delta$ in its
incident tetrahedron of $N_i$. By
Lemma~\ref{lem:one-boundary-triangle}, all edges $x_iy_j$ are
interior, while the edges of $\Delta$ remain boundary edges.
The three triangles $x_1x_2y_j$ therefore have only interior
edges and join the three new tetrahedra in a $3$-cycle of
the singular set of $P_{N_1\sharp N_2}$.
Every other face of a new tetrahedron contains an edge of
$\Delta$, so this $3$-cycle forms a separate connected
component of the singular set.
The removed tetrahedra contain no triangle with three interior
edges, and the remaining triangles retain their type.
Thus the singular set of $P_{N_1\sharp N_2}$ is the disjoint
union of the singular sets of $P_{N_1}$ and $P_{N_2}$
and an isolated $3$-cycle.
Consequently, $P_{N_1\sharp N_2}$ is a contractible
closed simple polyhedron with noncellular intrinsic
stratification.

\end{remark}

Counting nuclei through their intrinsic polyhedra
could provide an alternative route to
Theorem~\ref{thm:equivalent-exponential-counting}.
However, the examples above show that such an approach
is less straightforward than one might hope, since
these polyhedra need not be special spines.

\section{Dual graphs and genus}
\label{section:GraphGenus}
We now study the consequences of the structural results
and gluing operation from the previous section for the
dual graphs of nuclei.
For a triangulated $3$-ball $N$, recall that $G_N$ denotes its dual graph:
the vertices of $G_N$ are the tetrahedra of $N$, and two vertices are joined
when the corresponding tetrahedra share an interior triangle. In other words,
$G_N$ is the $1$-skeleton of the dual complex $D_N$. Recall also
that, for a graph $G$, $\gamma(G)$ denotes its orientable genus, that is, the
minimum genus of an orientable surface in which $G$ embeds. We write
\[
\gamma(N):=\gamma(G_N).
\]

\subsection{Nonplanarity of dual graphs}
Our earlier Theorem~\ref{thm:all-interior-tetrahedron} gives a direct
obstruction to planarity of the dual graph of a nucleus. Indeed, the fully interior
tetrahedron provided by the theorem determines a subdivision of $K_5$ in
$G_N$. Recall that a \emph{subdivision} of a graph is obtained by replacing
edges by paths with distinct new internal vertices.

\begin{theorem}
\label{thm:dual-graph-nonplanar}
Let $N$ be a nontrivial nucleus. Then $G_N$ contains a subdivision of $K_5$.
In particular,
\[
\gamma(G_N)\ge 1.
\]
\end{theorem}

\begin{proof}
    By Theorem~\ref{thm:all-interior-tetrahedron}, $N$ contains at least one tetrahedron with all six edges interior, call it $\sigma$. In particular, all
    four triangles of $\sigma$ are interior. Denote
    \begin{compactitem}
        \item \(    \tau_1,\tau_2,\tau_3,\tau_4\) the triangles of $\sigma$,
        \item $v$ the vertex of $G_N$ corresponding to $\sigma$,
        \item $v_i$ the vertex corresponding to the tetrahedron on the other side of $\tau_i$.
    \end{compactitem}
Thus $v$ is adjacent to each of $v_1,\ldots,v_4$.
    
    For every pair $i\ne j$, the triangles $\tau_i$ and $\tau_j$ meet in an
    interior edge $e_{ij}$ of $\sigma$. The tetrahedra containing $e_{ij}$ occur
    cyclically around this edge, and hence determine a cycle in $G_N$. This cycle
    contains $v$, with $v_i$ and $v_j$ as its two neighbors. Removing $v$
    from the cycle therefore leaves a path $P_{ij}$ from $v_i$ to $v_j$. Moreover, for $k\notin\{i,j\}$, the triangle $\tau_k$ and
    the edge $e_{ij}$ together contain all four vertices
    of $\sigma$. Thus, no tetrahedron other than $\sigma$
    contains both, and $v_k$ cannot lie on $P_{ij}$.
    
    We claim that these six paths are internally disjoint. Indeed, suppose that a tetrahedron
    other than $\sigma$ lies on the cycles corresponding to two distinct edges of
    $\sigma$. If the two edges meet, they span a triangle of $\sigma$, and the only
    tetrahedron other than $\sigma$ containing that triangle is the tetrahedron on
    the other side of it. This is precisely the common endpoint of the two
    corresponding paths. If the two edges are opposite, any tetrahedron containing
    both would have the same four vertices as $\sigma$, and hence would equal
    $\sigma$.

    Thus the subgraph consisting of the vertices
    \[
    v,v_1,v_2,v_3,v_4,
    \]
    the four edges $vv_i$, and the six paths $P_{ij}$ is a subdivision of
    $K_5$. Indeed, the five displayed vertices correspond to the five vertices
    of $K_5$; the edges from $v$ to the four vertices $v_i$ are already present
    in $G_N$, while for each pair $i\ne j$, the path $P_{ij}$ connects $v_i$ to
    $v_j$. By the preceding argument, these six paths have no interior vertices
    in common, and they meet the other four edges only at their endpoints. Hence
    $G_N$ contains a subdivision of $K_5$, and is therefore nonplanar.
    Consequently,
    \[
    \gamma(G_N)\ge 1.
    \]
\end{proof}

\subsection{Unbounded genus of nuclei}


\begin{figure}[tb]
    \centering
    \includegraphics[height=9em]{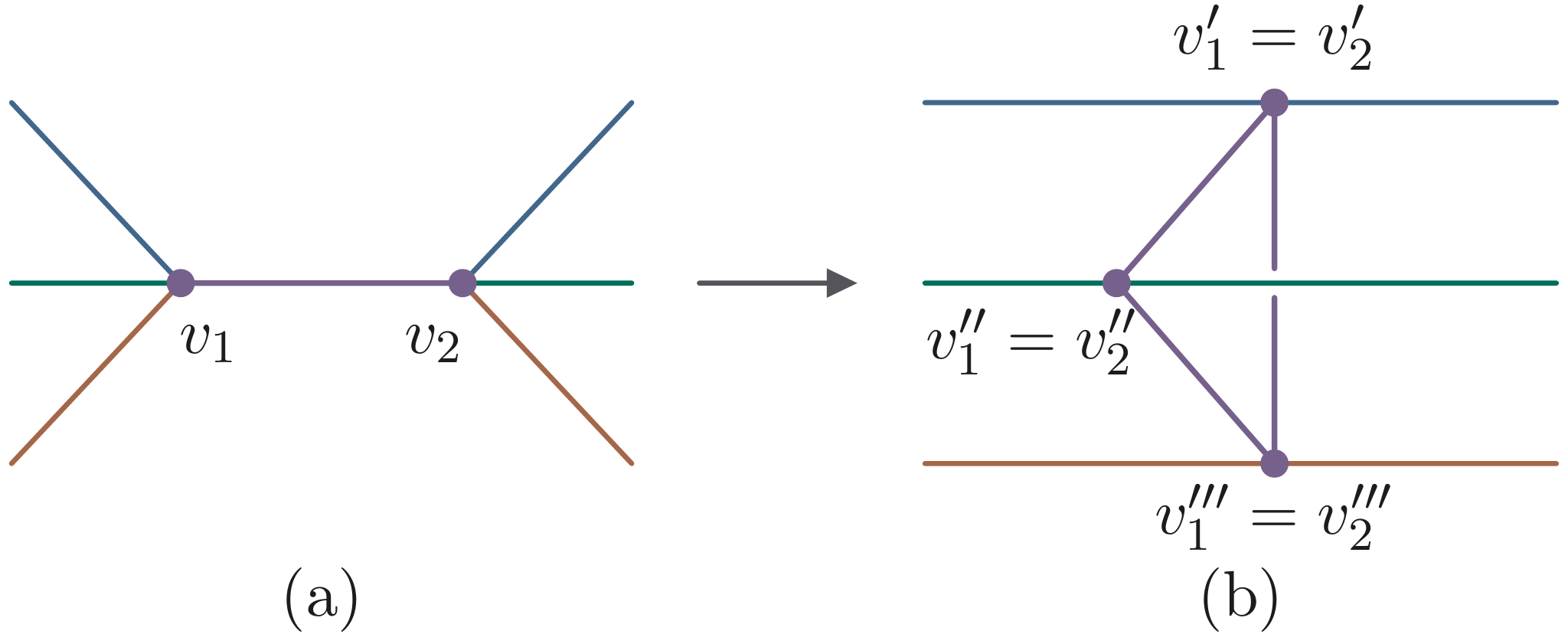}
    \caption{The local effect of the gluing operation on the dual graph.
    (a) After identifying the two boundary triangles and before the
    $2$-to-$3$ move. (b) After the $2$-to-$3$ move.}
    \label{fig:dual-gluing}
\end{figure}

We now use the gluing operation introduced in Section \ref{subsec:gluing-nuclei}. Recall that, given nuclei
$N_1$ and $N_2$ and chosen boundary triangles, $N_1\sharp N_2$ is obtained by
identifying the two triangles and performing the $2$-to-$3$ bistellar flip on
the resulting interior triangle. By Theorem~\ref{thm:gluenuclei}, if at least
one of $N_1,N_2$ is nontrivial, then $N_1\sharp N_2$ is again a nucleus, and
\[
T(N_1\sharp N_2)=T(N_1)+T(N_2)+1,
\]
where $T(N)$ denotes the number of tetrahedra of $N$.

To compare the genera of $N_1$, $N_2$, and $N_1\sharp N_2$, we now describe
how the gluing operation changes their dual graphs. Let $v_i\in G_{N_i}$ be the vertex corresponding to the
tetrahedron containing the boundary triangle chosen for gluing. For nontrivial $N_i$,
the vertex $v_i$ has degree $3$. Before the bistellar move, identifying the
boundary triangles adds an edge between $v_1$ and $v_2$, as shown on the left
of Figure~\ref{fig:dual-gluing}. After the move, the two vertices are replaced
by three vertices forming a triangle, as shown on the right of Figure~\ref{fig:dual-gluing}.

\begin{figure}[tb]
    \centering
    \includegraphics[height=9em]{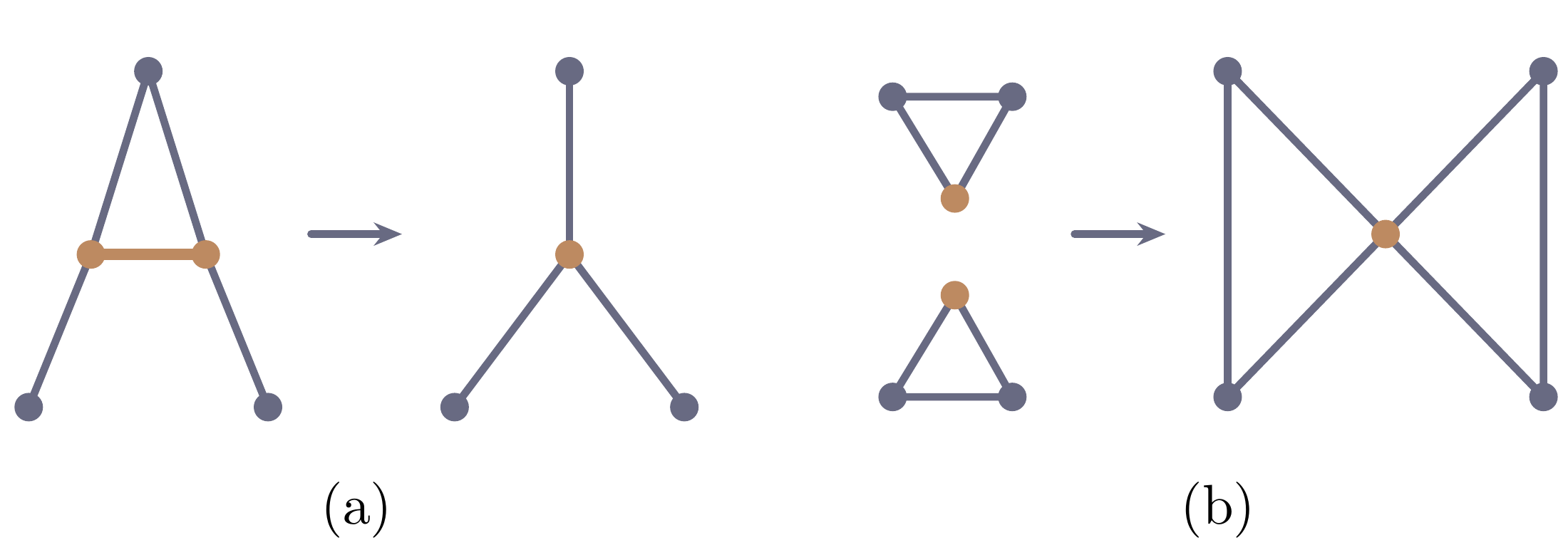}
    \caption{(a) Edge contraction, with the contracted edge highlighted; (b) vertex amalgamation at the highlighted vertices.
}
    \label{fig:graph_operations_pair}
\end{figure}

Moreover, we use \emph{graph minors}, obtained by deleting vertices or edges and \emph{contracting} edges. Contracting an edge means identifying its endpoints and removing the edge, discarding any resulting loops and duplicate edges. We also use the \emph{vertex amalgamation} $G_1\vee G_2$ of two disjoint graphs (see Figure~\ref{fig:graph_operations_pair}), obtained by identifying a chosen vertex of $G_1$ with a chosen vertex of $G_2$. The key observation here is that $G_{N_1}\vee G_{N_2}$ is a minor of $G_{N_1\sharp N_2}$.

\begin{proposition}
\label{prop:genus-superadditivity}
It holds for all nuclei $N_1$ and $N_2$ that
\[
\gamma(N_1\sharp N_2)\ge \gamma(N_1)+\gamma(N_2).
\]
\end{proposition}

\begin{proof}
In $G_{N_1\sharp N_2}$, the three new tetrahedra created by the
$2$-to-$3$ move correspond to the triangle in
Figure~\ref{fig:dual-gluing}. Delete one edge of this triangle and contract
the other two. The three vertices are then identified to a single vertex, and
the resulting graph is precisely the vertex amalgamation
$G_{N_1}\vee G_{N_2}$ obtained by identifying $v_1$ with $v_2$.

Thus $G_{N_1}\vee G_{N_2}$ is a minor of $G_{N_1\sharp N_2}$. Since
the orientable genus does not increase under taking minors \cite{BM2001}, and genus is additive
under vertex amalgamation \cite{BattleHararyKodamaYoungs},
\[
\gamma(N_1\sharp N_2)
\ge \gamma(G_{N_1}\vee G_{N_2})
=\gamma(G_{N_1})+\gamma(G_{N_2})
=\gamma(N_1)+\gamma(N_2).
\]
\end{proof}

This gives the following lower bound: the genus of nuclei can grow at least
linearly with the number of tetrahedra.

\begin{corollary}
\label{cor:linear-genus-growth}
There exist a constant $c>0$ and nuclei $N_k$, $k\ge1$, with
$T(N_k)\to\infty$ such that
\[
\gamma(N_k)\ge c\,T(N_k)
\]
for every $k$. 
\end{corollary}

\begin{proof}
Fix a nontrivial nucleus $M$ and write $q=T(M)$. By
Theorem~\ref{thm:dual-graph-nonplanar}, $\gamma(M)\ge1$. Let $N_k$ be obtained
by gluing $k$ copies of $M$ successively. Repeated application of
Proposition~\ref{prop:genus-superadditivity} gives
\[
\gamma(N_k)\ge k\gamma(M)\ge k,
\]
while
\[
T(N_k)=kq+(k-1)=k(q+1)-1.
\]
Therefore
\[
\gamma(N_k)\ge \frac{T(N_k)+1}{q+1}\ge \frac{T(N_k)}{q+1}.
\]
Taking $c=1/(q+1)$ proves the claim. 
\end{proof}

\begin{remark}
The linear growth in Corollary~\ref{cor:linear-genus-growth} is optimal up to the
constant. Indeed, if $N$ has $t$ tetrahedra, then $G_N$ has $t$ vertices and
at most $2t$ edges. A spanning tree of $G_N$ has $t-1$ edges, and each of the remaining
at most $t+1$ edges can be added using at most one handle. Hence
\[
\gamma(N)\le t+1.
\]
Thus, the maximal possible order of growth of the genus is linear in the number
of tetrahedra.
\end{remark}

\begin{corollary}\label{cor: unb1}
There exist nuclei with arbitrarily large genus.
\end{corollary}

\section{Examples of nuclei}
\label{section:Examples}
Small triangulated $3$-balls with particular combinatorial
and topological properties have attracted considerable
interest, particularly in connection with collapsibility,
evasiveness, local constructibility, and knotted edges
\cite{BL2013,BZ11}.
In this section, we investigate small nuclei exhibiting
such properties. We begin with the unique vertex-minimal
nontrivial nucleus, which has $11$ vertices and $35$
tetrahedra, and then examine two nuclei with $12$
vertices and $37$ tetrahedra.  The following local move produces further nuclei with the same
$f$-vector from a given one.

\begin{figure}[tb]
    \centering
    \includegraphics[height=9em]{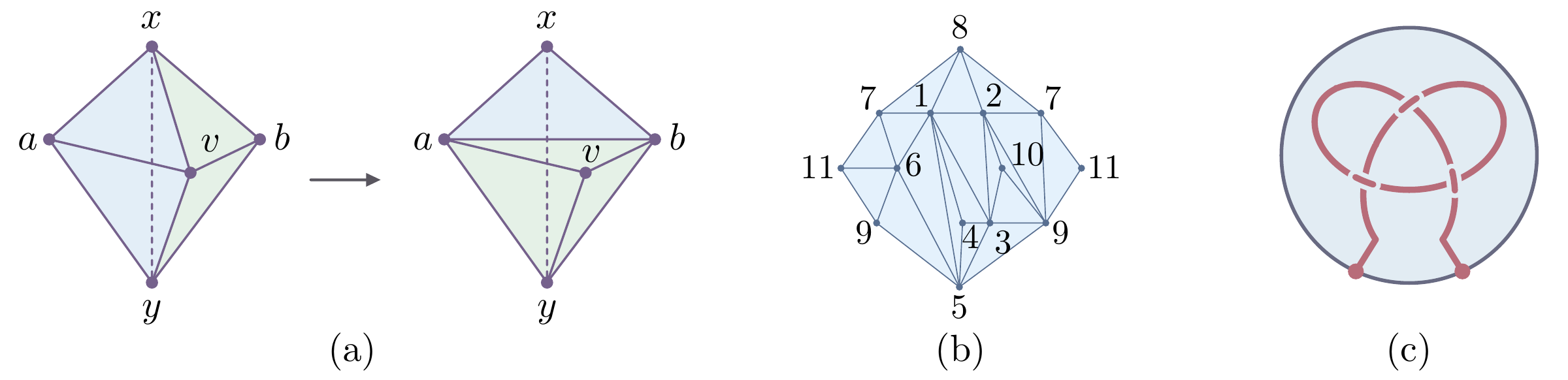}
    \caption{(a) The $2$-to-$2$ flip of tetrahedra;
    (b) The boundary triangulation of $N_{\min}$; (c) a knotted spanning edge.}
    \label{fig:nucleus-flip}
\end{figure}

Suppose that $[a,x,v]$ and $[b,x,v]$ are boundary triangles of a nucleus $N$
sharing the boundary edge $[x,v]$, and that their incident tetrahedra are
\[
[v,a,x,y]\qquad\text{and}\qquad[v,b,x,y].
\]
Thus the two tetrahedra share the interior triangle $[v,x,y]$.  If $[a,b]$ is
not an edge of $N$, replace these tetrahedra by
\[
[a,b,x,y]\qquad\text{and}\qquad[a,b,v,y].
\]
This is a $2$-to-$2$ flip of tetrahedra, inducing the
standard bistellar flip on the boundary \cite{Pachner1991}.
We call it simply a \emph{flip}; see Figure~\ref{fig:nucleus-flip}.

\begin{lemma}\label{lem:nucleus-flip}
A flip of a nucleus is again a nucleus.  Moreover, a flip preserves the
$f$-vector.
\end{lemma}

\begin{proof}
On the boundary this is the usual $2$-to-$2$ bistellar flip, replacing
$[a,x,v]$ and $[b,x,v]$ by $[a,b,x]$ and $[a,b,v]$. Thus the result is
again a triangulated $3$-ball with all vertices on the boundary, and the
$f$-vector is unchanged.

The old interior triangle $[v,x,y]$ contains the boundary edge $[x,v]$,
so $[x,y]$ and $[v,y]$ are interior edges.  Hence $[a,x,y]$ and $[b,x,y]$
are interior triangles.  Since $[a,x]$ and $[b,x]$ are boundary edges,
the nucleus condition forces $[a,y]$ and $[b,y]$ to be interior.
Thus the new interior triangle $[a,b,y]$ contains only the boundary edge
$[a,b]$, and all other interior triangles are unchanged.
\end{proof}

Several properties of our examples are verified by computer; the
scripts and data are in the reproducibility repository
\cite{ChatGPT2026Examples}. The README of each directory states the
topological arguments that turn the computations into proofs.

\subsection{The smallest nucleus}
\label{subsec:smallest-nucleus}
\begin{table}[t]
    \centering
    \small
    \setlength{\tabcolsep}{5pt}
    \begin{tabular}{|c|c|c|c|c|}
    \hline
    $\mathbf{1}\enspace \mathbf{2}\enspace \mathbf{3}\enspace 11$ & $\mathbf{1}\enspace \mathbf{7}\enspace \mathbf{8}\enspace 10$ & $2\enspace 4\enspace 6\enspace 11$ & $3\enspace 5\enspace 7\enspace 8$ & $4\enspace 6\enspace 8\enspace 10$ \\
    $\mathbf{1}\enspace \mathbf{2}\enspace \mathbf{8}\enspace 11$ & $1\enspace 8\enspace 10\enspace 11$ & $\mathbf{2}\enspace 4\enspace \mathbf{7}\enspace \mathbf{9}$ & $\mathbf{3}\enspace \mathbf{5}\enspace 8\enspace \mathbf{9}$ & $4\enspace \mathbf{6}\enspace \mathbf{9}\enspace \mathbf{11}$ \\
    $\mathbf{1}\enspace \mathbf{3}\enspace \mathbf{4}\enspace 11$ & $\mathbf{2}\enspace \mathbf{3}\enspace 6\enspace \mathbf{10}$ & $\mathbf{2}\enspace 4\enspace \mathbf{9}\enspace \mathbf{10}$ & $3\enspace 6\enspace 7\enspace 10$ & $4\enspace \mathbf{7}\enspace \mathbf{9}\enspace \mathbf{11}$ \\
    $\mathbf{1}\enspace \mathbf{4}\enspace \mathbf{5}\enspace 11$ & $2\enspace 3\enspace 6\enspace 11$ & $\mathbf{2}\enspace 5\enspace \mathbf{7}\enspace \mathbf{8}$ & $3\enspace \mathbf{6}\enspace \mathbf{7}\enspace \mathbf{11}$ & $4\enspace 8\enspace 9\enspace 10$ \\
    $\mathbf{1}\enspace \mathbf{5}\enspace \mathbf{6}\enspace 10$ & $2\enspace 4\enspace 5\enspace 7$ & $2\enspace 5\enspace 8\enspace 11$ & $3\enspace 7\enspace 8\enspace 10$ & $\mathbf{5}\enspace \mathbf{6}\enspace 8\enspace \mathbf{9}$ \\
    $1\enspace 5\enspace 10\enspace 11$ & $2\enspace 4\enspace 5\enspace 11$ & $\mathbf{3}\enspace \mathbf{4}\enspace \mathbf{5}\enspace 7$ & $\mathbf{3}\enspace 8\enspace \mathbf{9}\enspace \mathbf{10}$ & $5\enspace 6\enspace 8\enspace 10$ \\
    $\mathbf{1}\enspace \mathbf{6}\enspace \mathbf{7}\enspace 10$ & $2\enspace 4\enspace 6\enspace 10$ & $3\enspace 4\enspace 7\enspace 11$ & $4\enspace 6\enspace 8\enspace 9$ & $5\enspace 8\enspace 10\enspace 11$ \\
    \hline
    \end{tabular}
    \caption{The $35$ tetrahedra of $N_{\min}$. The vertices of each boundary triangle are shown in boldface.}
    \label{tab:Nmin-facets}
\end{table}
Let $N_{\min}$ be the simplicial complex on the vertex set
$\{1,\ldots,11\}$ whose $35$ tetrahedra are listed in Table~\ref{tab:Nmin-facets}.
The $f$-vectors of the triangulation and its boundary are \(f(N_{\min})=(11,54,79,35)\), and \(
f(\partial N_{\min})=(11,27,18)\).
The triangulation of the boundary is shown in Figure~\ref{fig:nucleus-flip}(b). The next theorem explains the notation $N_{\min}$.

\begin{theorem}
\label{thm:smallest-nucleus}
Every nontrivial nucleus has at least $11$ vertices. Up to simplicial
isomorphism, $N_{\min}$ is the unique nontrivial nucleus on $11$ vertices.
\end{theorem}

\begin{proof}
Brinkmann and McKay's \texttt{plantri}~\cite{BrinkmannMcKay2007}
enumerates triangulations of the $2$-sphere up to isomorphism.
We use it to enumerate the possible boundaries of nuclei.
There are $306$ such triangulations with $4$ to $10$ vertices,
and $1249$ with $11$ vertices.

For each boundary triangulation $S$, we consider all possible
sets of tetrahedra on its vertex set. Every boundary triangle
must belong to exactly one tetrahedron, and every other
triangle to either zero or two. We impose the nucleus
condition and the Euler characteristic condition for vertex
links. Moreover, by Theorem~\ref{thm:all-interior-tetrahedron},
at least one tetrahedron must have all six edges interior.

Since all vertices of $N$ lie on its boundary,
which is a triangulated $2$-sphere, we have
\[
f_1(\partial N)=3v-6,
\qquad
f_2(\partial N)=2v-4,
\]
where $v$ is the number of vertices of $N$.
Let $t$ be the number of its tetrahedra.
Counting triangles gives
\(
4t=2f_2(N)-f_2(\partial N),
\)
hence $f_2(N)=2t+v-2$. Euler's formula for the
$3$-ball then yields
\(
v-f_1(N)+f_2(N)-t=1,
\)
so $f_1(N)=t+2v-3$. Consequently,
\[
f_1^{\mathrm{int}}(N)
=f_1(N)-f_1(\partial N)
=t-v+3,
\]
where $f_1^{\mathrm{int}}(N)$ denotes the number of
interior edges of $N$.

Since $N$ is simplicial, $f_1(N)\le\binom{v}{2}$,
and therefore
\[
t\le\binom{v}{2}-2v+3.
\]
Moreover, by Theorem~\ref{thm:all-interior-tetrahedron},
$f_1^{\mathrm{int}}(N)\ge6$, so $t\ge v+3$.
Thus, for $v=11$, it suffices to consider
$14\le t\le36$.
The exhaustive computation finds no candidates with fewer
than $11$ vertices. For $v=11$, the computation yields $1829$ candidates,
of which only one is a triangulated $3$-ball,
namely $N_{\min}$.
The complete enumeration and verification are provided
in \texttt{classification/} of the reproducibility
archive~\cite{ChatGPT2026Examples}.
\end{proof}

\begin{figure}[tb]
    \centering
    \includegraphics[height=9em]{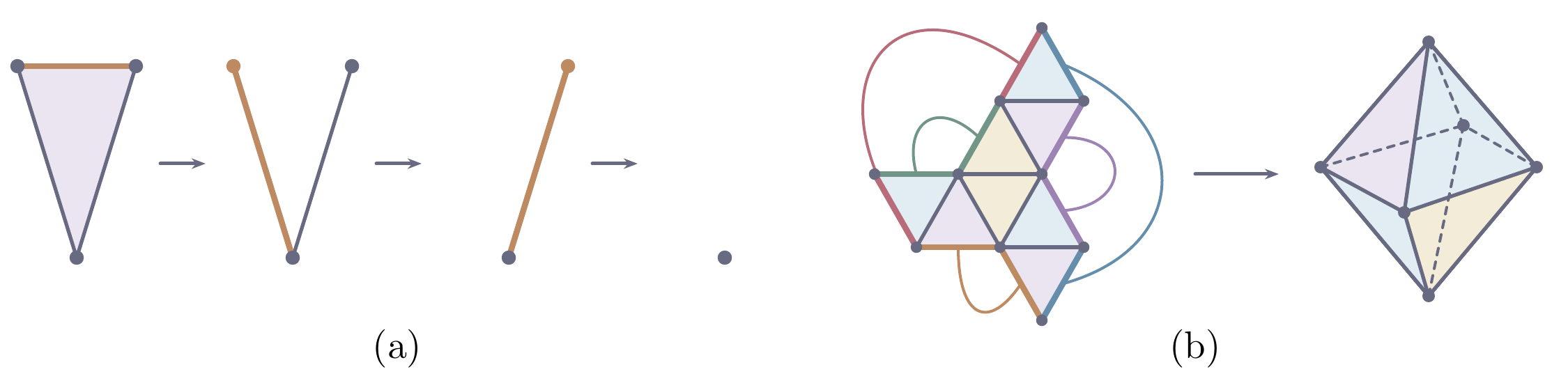}
    \caption{(a) A sequence of elementary simplicial collapses;
    (b) a local construction producing a triangulated $2$-sphere.}
    \label{fig:topology-illustrations}
\end{figure} 

In particular, $N_{\min}$ admits no flip: no boundary edge of $N_{\min}$ lies in only
two tetrahedra whose remaining vertices $a,b$ are not joined by an edge. Besides its minimality, $N_{\min}$ has several interesting
combinatorial and topological properties, summarized in the
following proposition. Let us recall a few standard definitions of notions that we will need in this result. 
Following \cite{BL2013}, a \emph{spanning edge} of a triangulated
$3$-ball is an edge with both endpoints on the boundary and
its relative interior in the interior of the ball. In a nucleus,
these are precisely its interior edges. Closing a spanning edge
with an arc in the boundary gives a knot, whose type is
independent of the chosen boundary arc. A spanning edge is
\emph{knotted} if this knot is nontrivial; see Figure~\ref{fig:nucleus-flip}(c).

Moreover, a complex is
\emph{collapsible} if it collapses to a vertex (see Figure~\ref{fig:topology-illustrations}(a)), and \emph{contractible}
if its realization is homotopy equivalent to a point
\cite{Hatcher2002}. Thus, collapsibility implies contractibility, but not conversely.

A simplicial complex is \emph{non-evasive}
if it is a nonempty
simplex, or if there is a vertex $v$ for which both its deletion
$\operatorname{del}_K(v)=\{\sigma\in K:v\notin\sigma\}$ and its link
$\operatorname{lk}_K(v)$ are non-evasive. Otherwise it is \emph{evasive}.
The recursion starts with a simplex containing at least one vertex;
a complex with no vertices is considered evasive.
A \emph{tree of tetrahedra} is a
$3$-ball obtained by successively attaching tetrahedra along single
boundary triangles; equivalently, its dual graph is a tree.
A $3$-ball is \emph{locally constructible}
(\emph{LC}) if it can be obtained from a tree of tetrahedra by successively
identifying pairs of boundary triangles that share a boundary edge,
with a simplicial $3$-ball as the final result
\cite{BenedettiZiegler2011LC}.
Figure~\ref{fig:topology-illustrations}(b) illustrates this
construction for a two-dimensional example.

\begin{figure}[tb]
    \centering
    \includegraphics[height=9em]{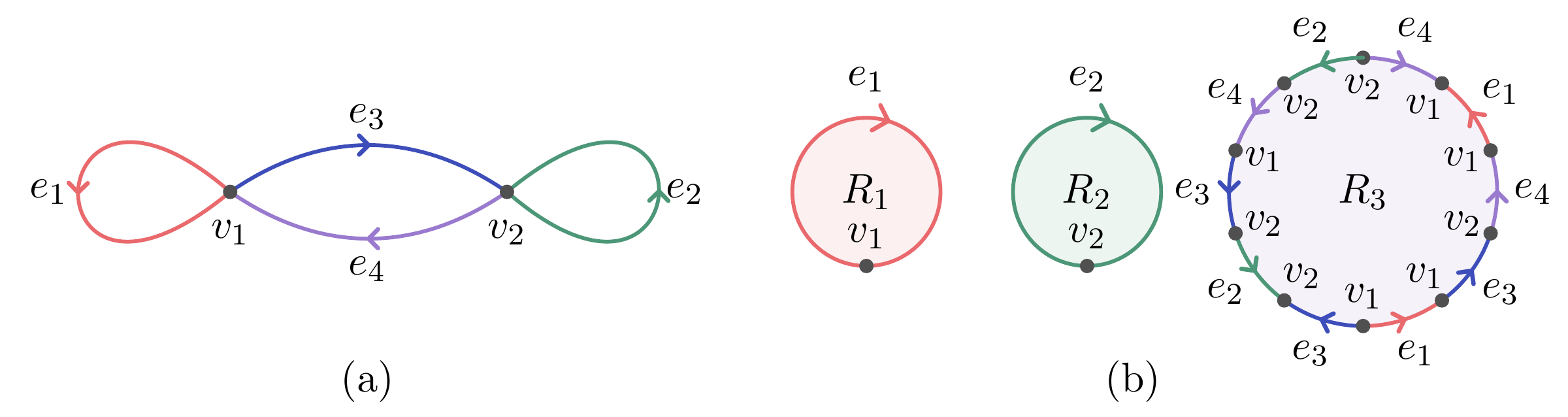}
    \caption{The intrinsic polyhedron $P_{N_{\min}}$:
        (a) its $1$-skeleton and (b) the attaching maps
        of its three regions.}
    \label{fig:DnPn}
\end{figure}

\begin{proposition}\label{prop:Nmin-properties}
The nucleus $N_{\min}$ has the following properties.
\begin{compactenum}
    \item It is collapsible and evasive, but not locally
      constructible.
    \item All its $27$ spanning edges are unknotted.
    \item Its dual graph has orientable genus $3$.
    \item Exactly two tetrahedra have all six edges interior,
          namely $[4,6,8,10]$ and $[5,8,10,11]$.
    \item Its dual complex $D_{N_{\min}}$ has $35$      vertices,
      $61$ edges, and $27$ two-cells, while its intrinsic
      polyhedron $P_{N_{\min}}$ has $2$ vertices,
      $4$ edges, and $3$ regions.  
    \item Its intrinsic polyhedron is cellular and homeomorphic     to Bing's house with two rooms \cite{FaganQiuWang2026}.     The polyhedron $P_{N_{\min}}$ is illustrated in
            Figure~\ref{fig:DnPn}.
\end{compactenum}
\end{proposition}

\begin{proof}
An explicit sequence of elementary collapses proving collapsibility of $N_{\min}$
is provided in repository \texttt{classification/} of the reproducibility
repository \cite{ChatGPT2026Examples}. Its evasiveness is verified computationally in
\texttt{nuclei\_intrinsic\_spines\_and\_evasiveness/} of the
same repository. By Benedetti~\cite[Corollary~36]{BB2017},
the tetrahedron is the only Mogami nucleus.
Since every locally constructible nucleus is Mogami,
$N_{\min}$ is not locally constructible.

The computations in \texttt{spanning-edges/} verify that all $27$
spanning edges are unknotted. The exact genus computation in
\texttt{dual-genus/} exhibits a genus-$3$ embedding of the dual graph
and excludes embeddings of lower genus.

A direct check shows that precisely the tetrahedra
$[4,6,8,10]$ and $[5,8,10,11]$ have all six edges interior.
The dual complex $D_{N_{\min}}$ has one vertex for each of the
$35$ tetrahedra, one edge for each of the $61$ interior triangles,
and one $2$-cell for each of the $27$ interior edges.

By Proposition~\ref{prop:intrinsic-structure}, the two
tetrahedra above correspond to the intrinsic vertices of
$P_{N_{\min}}$, whose $1$-skeleton is $4$-regular.
Counting loops twice, the degree-sum formula
\cite[Proposition~1.3.3]{West2001} gives $2e=4\cdot 2$,
hence $e=4$. Moreover, $P_{N_{\min}}$ is contractible
by Theorem~\ref{thm:dual-closed-simple}, and its intrinsic
stratification is cellular by
Proposition~\ref{prop:intrinsic-structure}.
Thus, Euler's formula gives $1=2-4+r$, and hence $r=3$ regions.

Comparing the attaching maps in Figure~\ref{fig:DnPn} with the classification of acyclic cellular fake surfaces of complexity two by Fagan, Qiu, and Wang~\cite{FaganQiuWang2026} identifies $P_{N_{\min}}$ as Bing's house with two rooms. The identification is verified computationally in \texttt{nuclei\_intrinsic\_spines\_and\_evasiveness/} of the repository \cite{ChatGPT2026Examples}.
\end{proof}

Starting from $N_{\min}$, we construct nine further  pairwise non-isomorphic nuclei with $12$ vertices and $37$ tetrahedra. Attach a tetrahedron $[a,b,c,y]$, with a new vertex $y$, along a boundary triangle $[a,b,c]$ of $N_{\min}$, whose incident tetrahedron is $[a,b,c,x]$. The \emph{$2$-to-$3$ bistellar move} (see Theorem \ref{thm:gluenuclei}) replaces these two tetrahedra by
\[
[a,b,x,y],\qquad [a,c,x,y],\qquad [b,c,x,y],
\]
provided $[x,y]$ is not already an edge, as is automatically the case here.

The nine choices arise from a symmetry of $N_{\min}$ that pairs up its $18$ boundary triangles, given by the permutation
\[
(1\ 9)(3\ 7)(4\ 11)(5\ 6)(8\ 10),
\]
with vertex $2$ fixed. All spanning edges of the resulting nuclei are unknotted, and the two tetrahedra of $N_{\min}$ with all six edges interior remain unchanged. The constructions and their verification are available in \texttt{spanning-edges/} of \cite{ChatGPT2026Examples}.

\subsection{Other small examples}
\label{subsec:other-small-examples}

Before the discovery of $N_{\min}$, the smallest known
nontrivial nuclei had $12$ vertices and $37$ tetrahedra.
We present two such examples with the same $f$-vector as the nine nuclei above, but with different structural properties.

\begin{table}[t]
    \centering
    \small
    \setlength{\tabcolsep}{5pt}
    \begin{tabular}{|c|c|c|c|c|}
    \hline
    $\mathbf{1}\enspace \mathbf{3}\enspace \mathbf{4}\enspace 10$ & $1\enspace 6\enspace 7\enspace 10$ & $\mathbf{2}\enspace \mathbf{7}\enspace 11\enspace \mathbf{12}$ & $\mathbf{4}\enspace \mathbf{5}\enspace \mathbf{6}\enspace 11$ & $4\enspace 8\enspace 10\enspace 11$ \\
    $1\enspace 3\enspace 5\enspace 10$ & $\mathbf{2}\enspace 3\enspace \mathbf{5}\enspace \mathbf{9}$ & $\mathbf{2}\enspace 8\enspace \mathbf{9}\enspace \mathbf{10}$ & $\mathbf{4}\enspace \mathbf{5}\enspace \mathbf{7}\enspace 8$ & $6\enspace 7\enspace 8\enspace 9$ \\
    $\mathbf{1}\enspace \mathbf{3}\enspace 5\enspace \mathbf{11}$ & $2\enspace 3\enspace 5\enspace 11$ & $2\enspace 8\enspace 10\enspace 11$ & $4\enspace 5\enspace 8\enspace 11$ & $6\enspace 7\enspace 9\enspace 11$ \\
    $\mathbf{1}\enspace \mathbf{4}\enspace \mathbf{6}\enspace 10$ & $2\enspace 3\enspace 8\enspace 9$ & $3\enspace 4\enspace 9\enspace 10$ & $4\enspace 6\enspace 10\enspace 11$ & $6\enspace 7\enspace 10\enspace 11$ \\
    $1\enspace 5\enspace 7\enspace 8$ & $2\enspace \mathbf{3}\enspace \mathbf{8}\enspace \mathbf{11}$ & $\mathbf{3}\enspace \mathbf{4}\enspace 9\enspace \mathbf{12}$ & $4\enspace 7\enspace 8\enspace 9$ & $\mathbf{6}\enspace \mathbf{8}\enspace 9\enspace \mathbf{12}$ \\
    $1\enspace \mathbf{5}\enspace \mathbf{7}\enspace \mathbf{10}$ & $\mathbf{2}\enspace \mathbf{5}\enspace \mathbf{6}\enspace 11$ & $3\enspace \mathbf{5}\enspace \mathbf{9}\enspace \mathbf{10}$ & $\mathbf{4}\enspace \mathbf{7}\enspace 9\enspace \mathbf{12}$ & $6\enspace 9\enspace 11\enspace 12$ \\
    $\mathbf{1}\enspace 5\enspace \mathbf{8}\enspace \mathbf{11}$ & $\mathbf{2}\enspace \mathbf{6}\enspace 11\enspace \mathbf{12}$ & $\mathbf{3}\enspace \mathbf{8}\enspace 9\enspace \mathbf{12}$ & $4\enspace 8\enspace 9\enspace 10$ & $7\enspace 9\enspace 11\enspace 12$ \\
    $\mathbf{1}\enspace \mathbf{6}\enspace 7\enspace \mathbf{8}$ & $\mathbf{2}\enspace \mathbf{7}\enspace \mathbf{10}\enspace 11$ &  &  &  \\
    \hline
    \end{tabular}
    \caption{The $37$ tetrahedra of $N_{\mathrm{CEY}}$. The vertices of each boundary triangle are shown in boldface.}
    \label{tab:NCEY-facets}
\end{table}

\begin{figure}[tb]
    \centering
    \begin{subfigure}[b]{0.3\textwidth}
        \centering
        \includegraphics[height=9em]{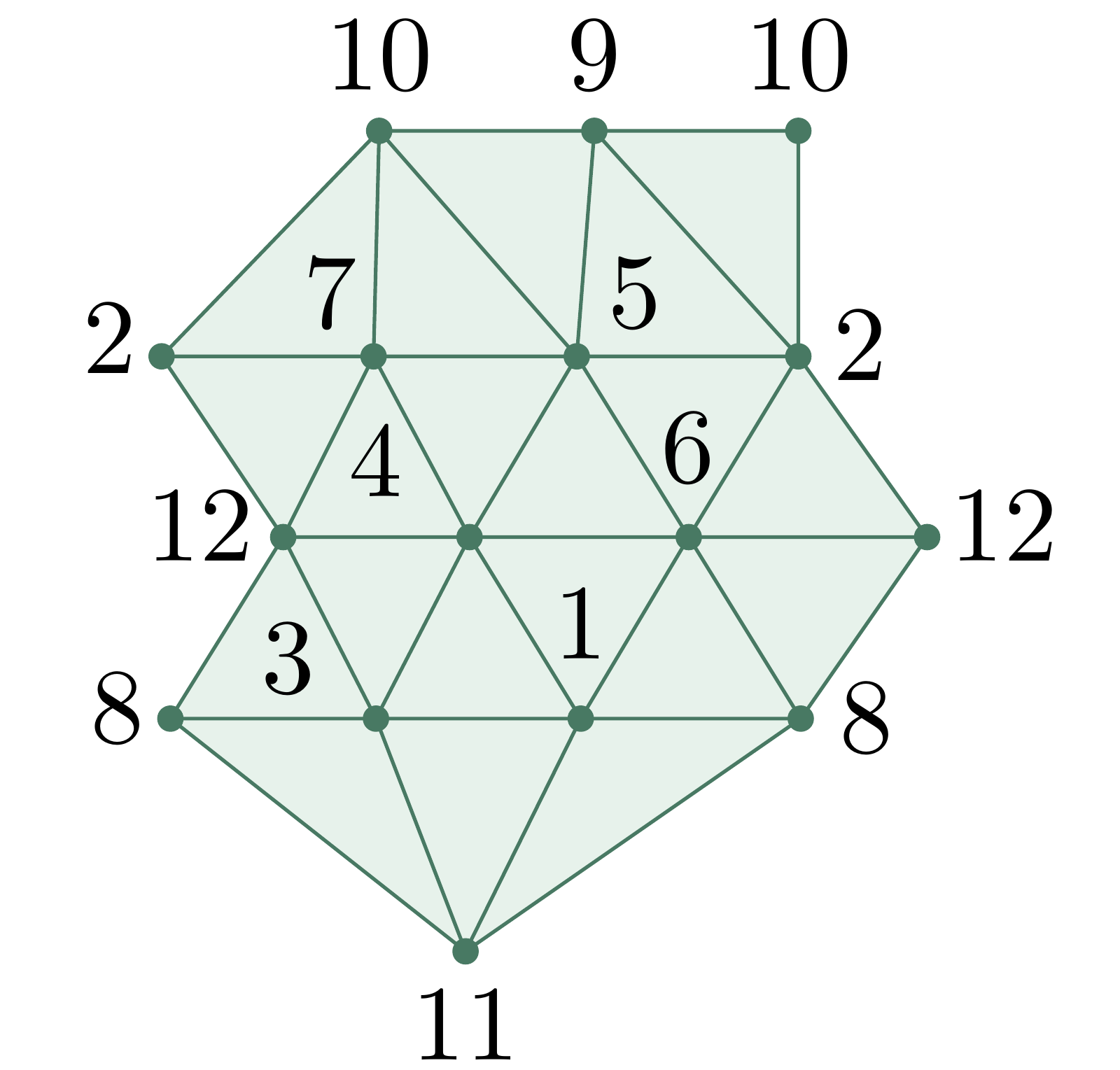}
        \caption{}
        \label{fig:first}
    \end{subfigure}
    \hfill
    \begin{subfigure}[b]{0.3\textwidth}
        \centering
        \includegraphics[height=9em]{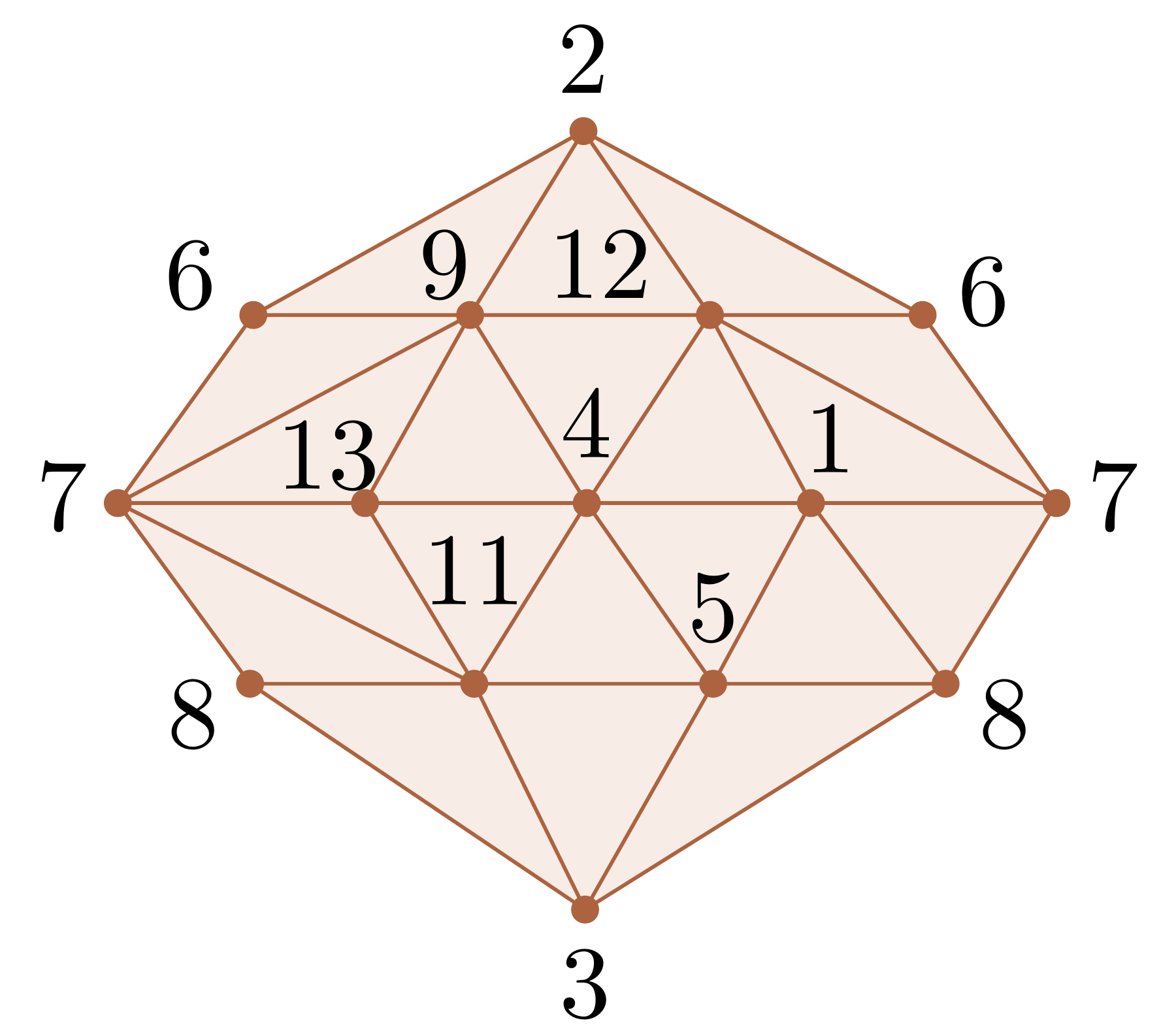}
        \caption{}
        \label{fig:second}
    \end{subfigure}
    \hfill
    \begin{subfigure}[b]{0.3\textwidth}
        \centering
        \includegraphics[height=5em]{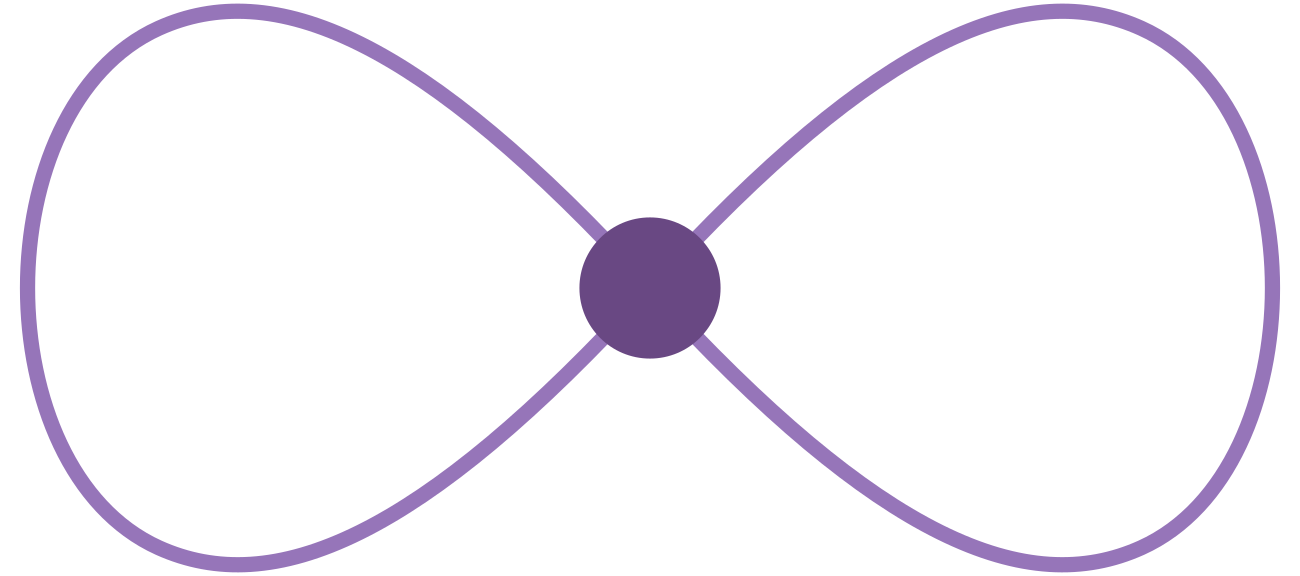}
        \caption{}
        \label{fig:third}
    \end{subfigure}
    \caption{(a) Boundary triangulation of the $12$-vertex nucleus $N_{\mathrm{CEY}}$, (b) and of $N_{\mathrm{BL}}$; (c) The abalone is the only special spine of $B^3$ with one
intrinsic vertex. Its 1-skeleton is depicted here.}
    \label{fig:both}
\end{figure}

\paragraph{The Collet--Eckmann--Younan nucleus.}
Collet, Eckmann, and Younan obtained their $12$-vertex nucleus $N_{\mathrm{CEY}}$
by applying their computer-assisted reduction procedure
to a large triangulated $3$-ball containing a trefoil knot
\cite{CEY2013}.  
The $f$-vectors of the nucleus and its boundary are \(f(N_{\mathrm{CEY}})=(12,58,84,37)\), and \(
f(\partial N_{\mathrm{CEY}})=(12,30,20)\). Its tetrahedra are reproduced in Table~\ref{tab:NCEY-facets},
and its boundary triangulation is shown in Figure~\ref{fig:both}(a). 

\paragraph{A nucleus from the Benedetti--Lutz ball.}
Let $N_{\mathrm{BL}}$ denote the nucleus obtained from the
non-LC ball $B_{12,38}$ of Benedetti and Lutz
\cite[Section~3]{BL2013} using the reduction moves of
\cite{CEY2013}.
The reduction begins by splitting vertex $10$ along a path,
in the terminology of \cite{CEY2013}.
In the resulting nucleus, the original interior vertex $4$
lies on the boundary, vertex $10$ has disappeared, and a new
vertex $1$ has been introduced. Its vertex set is therefore
$\{1,\ldots,9,11,12,13\}$.
Its tetrahedra are listed in Table~\ref{tab:NBL-facets},
and its boundary triangulation is shown in Figure~\ref{fig:both}(b).
\begin{table}[H]
    \centering
    \small
    \setlength{\tabcolsep}{5pt}
    \begin{tabular}{|c|c|c|c|c|}
    \hline
    $1\enspace 2\enspace 3\enspace 4$ & $2\enspace 4\enspace 5\enspace 7$ & $2\enspace 8\enspace 9\enspace 11$ & $\mathbf{3}\enspace \mathbf{5}\enspace \mathbf{8}\enspace 13$ & $3\enspace 8\enspace 12\enspace 13$ \\
    $1\enspace 2\enspace 3\enspace 7$ & $2\enspace 5\enspace 7\enspace 13$ & $\mathbf{2}\enspace 8\enspace \mathbf{9}\enspace \mathbf{12}$ & $\mathbf{3}\enspace \mathbf{5}\enspace 9\enspace \mathbf{11}$ & $4\enspace 5\enspace 6\enspace 7$ \\
    $\mathbf{1}\enspace 2\enspace \mathbf{4}\enspace \mathbf{5}$ & $2\enspace 5\enspace 8\enspace 13$ & $2\enspace 8\enspace 12\enspace 13$ & $3\enspace 5\enspace 9\enspace 13$ & $\mathbf{4}\enspace \mathbf{5}\enspace 6\enspace \mathbf{11}$ \\
    $\mathbf{1}\enspace 2\enspace \mathbf{5}\enspace \mathbf{8}$ & $\mathbf{2}\enspace \mathbf{6}\enspace \mathbf{9}\enspace 11$ & $3\enspace 4\enspace 6\enspace 7$ & $3\enspace \mathbf{6}\enspace \mathbf{7}\enspace \mathbf{12}$ & $\mathbf{4}\enspace 6\enspace \mathbf{11}\enspace \mathbf{13}$ \\
    $\mathbf{1}\enspace 2\enspace \mathbf{7}\enspace \mathbf{8}$ & $2\enspace 6\enspace 11\enspace 13$ & $3\enspace 4\enspace 6\enspace 13$ & $3\enspace 6\enspace 12\enspace 13$ & $5\enspace \mathbf{6}\enspace \mathbf{7}\enspace \mathbf{9}$ \\
    $\mathbf{1}\enspace 3\enspace \mathbf{4}\enspace \mathbf{12}$ & $\mathbf{2}\enspace \mathbf{6}\enspace \mathbf{12}\enspace 13$ & $3\enspace \mathbf{4}\enspace \mathbf{9}\enspace \mathbf{12}$ & $\mathbf{3}\enspace \mathbf{8}\enspace 9\enspace \mathbf{11}$ & $5\enspace 6\enspace 9\enspace 11$ \\
    $\mathbf{1}\enspace 3\enspace \mathbf{7}\enspace \mathbf{12}$ & $2\enspace \mathbf{7}\enspace \mathbf{8}\enspace \mathbf{11}$ & $3\enspace \mathbf{4}\enspace \mathbf{9}\enspace \mathbf{13}$ & $3\enspace 8\enspace 9\enspace 12$ & $5\enspace \mathbf{7}\enspace \mathbf{9}\enspace \mathbf{13}$ \\
    $2\enspace 3\enspace 4\enspace 7$ & $2\enspace \mathbf{7}\enspace \mathbf{11}\enspace \mathbf{13}$ &  &  &  \\
    \hline
    \end{tabular}
    \caption{The $37$ tetrahedra of $N_{\mathrm{BL}}$. The vertices of each boundary triangle are shown in boldface.}
    \label{tab:NBL-facets}
\end{table}

\begin{proposition}\label{prop:other-small-properties}
The nuclei $N_{\mathrm{CEY}}$ and $N_{\mathrm{BL}}$ have the
following properties.
\begin{compactenum}
    \item Both contain a trefoil-knotted spanning edge,
          namely $[9,11]$ and $[2,3]$, respectively.
    \item Both nuclei are evasive.
    \item Their dual complexes have the same $f$-vector
          $(37,64,28)$, but their dual graphs have
          orientable genera $3$ and $4$, respectively.
          In particular, the nuclei are not combinatorially
          isomorphic.
    \item Each contains exactly one tetrahedron with all
          six edges interior, namely $[6,7,9,11]$ and
          $[2,3,4,7]$, respectively.
    \item Their intrinsic polyhedra both have one vertex,
          two edges, and two regions, and are both homeomorphic to
      the abalone (Figure \ref{fig:both}(c)) \cite{FaganQiuWang2026}.
\end{compactenum}
\end{proposition}

\begin{proof}
The trefoil-knotted edges  are verified computationally in \texttt{trefoil-edges/} and the dual-graph genera in \texttt{dual-genus/} of the reproducibility repository
\cite{ChatGPT2026Examples}. Their evasiveness is verified computationally in
\texttt{nuclei\_intrinsic\_spines\_and\_evasiveness/} of
the same repository.

The remaining properties follow by direct inspection of the facet
lists. In particular, $P_{N_{\mathrm{CEY}}}$ and $P_{N_{\mathrm{BL}}}$
each have one intrinsic vertex and a connected singular graph. By
Corollary~\ref{cor:special-spine}, they are therefore special spines
of $B^3$, and the abalone is the only special spine of $B^3$ with one
intrinsic vertex \cite[Table~2.1]{FaganQiuWang2026}.
\end{proof}

Let us remark that unlike $N_{\min}$, these nuclei are not rigid: $N_{\mathrm{CEY}}$ and
$N_{\mathrm{BL}}$ admit $6$ and $4$ flips, respectively, each producing
another nucleus with the same $f$-vector, and combinations of the flips even more.

\begin{remark}\label{rem:nuclei-all-sizes}
There exists a nucleus with $t$ tetrahedra for every
$t\ge 37$. Indeed, in $N_{\mathrm{CEY}}$, the interior
triangle $[3,4,10]$ is shared by the tetrahedra
$[1,3,4,10]$ and $[3,4,9,10]$.
Since $[1,9]$ is not an edge, a $2$-to-$3$
bistellar move replaces them by
\[
[1,3,4,9],\quad [1,3,9,10],\quad
[1,4,9,10].
\]
The boundary is unchanged, and each new interior
triangle contains exactly one boundary edge.
Thus the resulting triangulation is a nucleus
with $38$ tetrahedra.
By Theorem~\ref{thm:gluenuclei}, gluing a
tetrahedron to a nontrivial nucleus increases
its number of tetrahedra by two. Starting from
the nuclei with $37$ and $38$ tetrahedra,
we obtain the claim.
\end{remark}

\begin{remark}\label{rem:atoms}
A further reduction of nuclei was studied by
Younan~\cite{MY2014phd}; see also~\cite{Eckmann2015Atoms}.
An \emph{atom} is a nucleus with no collapsible boundary edges,
that is, a nucleus $N$ such that
\[
\operatorname{lk}_N(uv)\neq
\operatorname{lk}_N(u)\cap\operatorname{lk}_N(v)
\]
for every boundary edge $uv$.
The smallest nucleus $N_{\min}$ and the two $12$-vertex
nuclei $N_{\mathrm{CEY}}$~\cite{CEY2013} and
$N_{\mathrm{BL}}$ are all atoms.

Younan~\cite[Lemma~4.1.4]{MY2014phd} states that every
nucleus reduces to a unique atom, independently of the
order of boundary-edge contractions. However, this
assertion is not valid in general.
Indeed, we construct a nucleus with $20$ vertices and
$74$ tetrahedra admitting two collapsible boundary edges
whose contractions yield nonisomorphic atoms with
$71$ and $72$ tetrahedra, respectively.
The example and its verification are provided
in~\cite{ChatGPT2026Examples}.
The difficulty in the original argument is that
contracting an edge can change the admissibility of
other contractions. This observation does not rule out
the proposed reduction of the exponential-counting
problem from nuclei to atoms, but its justification
requires reconsideration.

We have also verified that the gluing operation of
Section~\ref{subsec:gluing-nuclei} preserves atoms,
provided both factors are nontrivial.
Consequently, the genus results from Section \ref{section:GraphGenus} extend to atoms.
Moreover, repeatedly gluing copies of $N_{\min}$ to
$N_{\mathrm{CEY}}$, along boundary triangles away from
its trefoil-knotted spanning edge, produces infinitely
many atoms containing an edge of the same knot type.
This is noteworthy in view of Younan's suggestion
\cite{MY2014phd} that the number of atoms containing
an edge of a prescribed knot type might be very small,
perhaps even one.
\end{remark}

\section*{Acknowledgements}
Bergfinnur Durhuus and Marta Pavelka were supported by the VILLUM Foundation $\#10059$.

\subsection*{Use of AI tools}
The computer code and verification scripts in the reproducibility
repository~\cite{ChatGPT2026Examples} were written with the help of
ChatGPT (OpenAI). During the preparation of the manuscript, ChatGPT
and Claude (Anthropic) were also used to check proofs, citations, and to independently re-implement several
computations, including the classification of $11$-vertex nuclei.
The authors take full responsibility for the
content of this paper.

\newpage

\bibliography{ResearchBibliography}
\bibliographystyle{amsplain}

\end{document}